\documentclass[11pt]{article}
\usepackage{latexsym}
\usepackage[all]{xy}
\usepackage{amsmath}
\usepackage{amssymb}
\usepackage{amsfonts}

\newtheorem{thm}{Th\'eor\`eme}[section]
\newtheorem{prop}[thm]{Proposition}
\newtheorem{lem}[thm]{Lemme}
\newtheorem{sublem}[thm]{Sous-lemme}
\newtheorem{df}[thm]{D\'efinition}
\newtheorem{cor}[thm]{Corollaire}

\newtheorem{rmk}[thm]{Remarque}

\begin{document}

\title{\textbf{Finitude homotopique des 
dg-alg\`ebres propres et lisses}}
\bigskip
\bigskip

\author{\bigskip\\
Bertrand To\"en\\
\small{Laboratoire Emile Picard UMR CNRS 5580} \\
\small{Universit\'e Paul Sabatier Toulouse 3, Bat 1R2}\\
\small{31 062 TOULOUSE cedex 9, France }}

\date{Septembre 2006}

\maketitle

\begin{abstract}
On montre que toute dg-alg\`ebre propre et lisse (sur un anneau
de base $k$) est 
d\'etermin\'ee \`a quasi-isomorphisme pr\`es par 
sa $\mathcal{A}_{n}$-alg\`ebre sous-jacente pour un certain $n$. 
De m\^eme, tout morphisme entre dg-alg\`ebres propres et lisses 
est d\'etermin\'e \`a homotopie pr\`es par le morphisme induit
sur les $\mathcal{A}_{n}$-alg\`ebres sous-jacentes.
On d\'emontre de plus que si $k$ est local alors l'entier $n$ peut \^etre choisi uniformement pour
toutes les dg-alg\`ebres propres et lisses dont deux 
invariants num\'eriques (le \emph{type} et la \emph{dimension cohomologique}) 
sont born\'es.
\end{abstract}

\medskip

\tableofcontents

\bigskip

\section{Introduction}

Soit $k$ un anneau commutatif. Une $k$-dg-alg\`ebre (associative et unitaire) $B$ peut 
aussi \^etre consid\'er\'ee comme une $\mathcal{A}_{n}$-alg\`ebre (pour laquelle
$\mu_{i}=0$ pour $i>2$, voir par exemple \cite{ke1,lef,ma}).  L'objectif de cet article 
est d'\'etudier sous quelles conditions sur $B$ existe-t-il un entier $n$ tel que $B$ soit 
determin\'ee \`a quasi-isomorphisme pr\`es par sa $\mathcal{A}_{n}$-alg\`ebre
sous-jacente. Les r\'esultats principaux de ce travail affirmentque cela 
est toujours le cas lorsque la dg-alg\`ebre $B$ est propre et lisse, c'est \`a dire
lorsque $B$ est parfait comme complexe de $k$-modules et aussi comme
un bi-$B$-dg-module (voir \cite{ke2,koso,tova1}). Plus pr\'ecis\`emment, notre prermier r\'esultat est 
le th\'eor\`eme suivant.

\begin{thm}{(voir Cor.  \ref{c2} et Cor. \ref{c4})}\label{ti1}
Soit $B$ et $B'$ deux $k$-dg-alg\`ebres propres et lisses sur $k$, et pour un entier $m$ notons
$i_{m}^{*}(B)$ et $i_{m}^{*}(B')$ les $\mathcal{A}_{m}$-alg\`ebres sous-jacentes.
Alors il existe un entier $n$ v\'erifiant les conditions suivantes.
\begin{enumerate}
\item Si $i_{n}^{*}(B)$ et $i_{n}^{*}(B')$ sont quasi-isomorphes comme $\mathcal{A}_{n}$-alg\`ebres, 
alors $B$ et $B'$ sont quasi-isomorphes comme dg-alg\`ebres.
\item L'application 
$$[B,B'] \longrightarrow [i_{n}^{*}(B),i_{n}^{*}(B')]$$
de l'ensemble des morphismes dans la cat\'egorie homotopique des dg-alg\`ebres
vers l'ensemble des morphismes dans la cat\'egorie homotopique des $\mathcal{A}_{n}$-alg\`ebres
est injective. 
\end{enumerate}
\end{thm}

Les r\'esultats que nous d\'emontrons sont en r\'ealit\'e plus pr\'ecis que
l'\'enonc\'e pr\'ec\'edent, car ils affirment que les morphismes
naturels sur les espaces de morphimes et les espaces d'\'equivalences
$$Map(B,B') \longrightarrow Map(i_{n}^{*}(B),i_{n}^{*}(B'))$$
$$Map^{eq}(B,B') \longrightarrow Map^{eq}(i_{n}^{*}(B),i_{n}^{*}(B')),$$
poss\`ede des r\'etractions \`a homotopie pr\`es. 

Le th\'eor\`eme \ref{ti1} est une intepr\'etation du fait que les dg-alg\`ebres
propres et lisses sont determin\'ees par leurs
$\mathcal{A}_{n}$-alg\`ebres sous-jacentes (pour un $n$ qui varie avec 
les dg-alg\`ebres consid\'er\'ees). C'est pour cette propri\'et\'e remarquable
que nous utilisons l'expression \emph{finitude homotopique} qui apparait dans
notre titre. 

Le second r\'esultat de ce travail pr\'ecise que l'entier $n$
peut \^etre choisit uniform\'ement pour toutes les dg-alg\`ebres propres et lisses
dont deux invariants num\'eriques, que nous appelons le \emph{type} et la \emph{dimension cohomologique},
sont born\'es. Un type est une application \`a support fini $\nu : \mathbb{Z} \longrightarrow \mathbb{N}$.
Une $k$-dg-alg\`ebre $B$ qui est propre est alors \emph{de type $\nu$}
si pour tout morphisme $k \longrightarrow K$ avec $K$ un corps, on a pour tout $i\in \mathbb{Z}$ 
$$Dim_{K}H^{i}(B\otimes_{k}^{\mathbb{L}}K)\leq \nu(i).$$
Pour une $k$-dg-alg\`ebre lisse $B$, nous introduisons un second invariant $d\in \mathbb{N}$, appel\'e
la \emph{dimension cohomologique}, qui mesure d'une certaine fa\c{c}on la longueur 
d'une r\'esolution libre de $B$ comme un bi-$B$-dg-module (on renvoie \`a la d\'efinition \ref{d4} pour les
d\'etails sur cette notion). Notre second th\'eor\`eme s'\'ennonce alors comme suit.

\begin{thm}{(voir Cor. \ref{c3})}\label{ti2} 
Supposons que $k$ soit un anneau local. Alors, pour tout
type $\nu$ et tout entier $d$, il existe un entier $n(\nu,d)$ 
tel que pour toutes $k$-dg-alg\`ebres $B$ et $B'$ de 
type $\nu$ et de dimension cohomologique inf\'erieure \`a $d$ on ait 
les deux propri\'et\'es suivantes.
\begin{enumerate}
\item Si $i_{n}^{*}(B)$ et $i_{n}^{*}(B')$ sont quasi-isomorphes comme $\mathcal{A}_{n}$-alg\`ebres, 
alors $B$ et $B'$ sont quasi-isomorphes comme dg-alg\`ebres.
\item L'application 
$$[B,B'] \longrightarrow [i_{n}^{*}(B),i_{n}^{*}(B')]$$
est injective. 
\end{enumerate}
\end{thm}

Tout comme le th\'eor\`eme \ref{ti1} nous d\'emontrerons en r\'ealit\'e 
un \'enonc\'e plus pr\'ecis portant sur les espaces de morphismes
et les espaces d'\'equivalences. De plus, lorsque $k$ n'est plus local 
nous montrons que le th\'eor\`eme \ref{ti2} reste valable localement
pour la topologie de Zariski sur $Spec\, k$ (voir Thm. \ref{t3}). Enfin, la preuve que nous
donnons du th\'eor\`eme \ref{ti2} montre que l'entier
$n(\nu,d)$ ne d\'epend que du couple $(\nu,d)$ et non pas de l'anneau $k$ (voir Thm. \ref{t4}).

Avant de donner quelques id\'ees des preuves des th\'eor\`emes \ref{ti1} et \ref{ti2}
mentionons le corollaire suivant. Pour cela, on rappelle qu'une $\mathcal{A}_{\infty}$-alg\`ebre
$B$ sur un corps $k$ est minimale si la diff\'erentielle de son complexe sous-jacent est nulle (voir \cite{ke1,lef}). 
On montre alors (voir \cite{lef}) que deux $\mathcal{A}_{\infty}$-alg\`ebres minimales sont isomorphes dans
$Ho(k-\mathcal{A}_{\infty}-alg)$ si et seulement si elles sont $\mathcal{A}_{\infty}$-isomorphes
(i.e. isomorphes dans une certaine cat\'egorie de $\mathcal{A}_{\infty}$-alg\`ebres
et $\mathcal{A}_{\infty}$-morphismes). 
On d\'eduit alors du th\'eor\`eme \ref{ti2} le corollaire suivant, qui est
une autre incarnation plus concr\`ete de la propri\'et\'e de finitude homotopique.

\begin{cor}\label{ci}
Pour tout type $\nu$  et tout entier $d\geq 0$ il existe un entier 
$n(\nu,d)$ qui poss\`ede la propri\'et\'e suivante. Pour tout
corps $k$, si $\{\mu_{i}\}_{i\geq 2}$ et $\{\mu'_{i}\}_{i\geq 2}$ sont deux structures
de $\mathcal{A}_{\infty}$-alg\`ebres (sur k) sur un m\^eme complexe \`a diff\'erentielle nulle fix\'e $V$, 
telles que les $\mathcal{A}_{\infty}$-alg\`ebres correspondantes $B$ et $B'$ soient
propres et lisses, et si $\mu_{i}=\mu'_{i}$ pour tout $i\leq n(\nu,d)$, alors
$B$ et $B'$ sont $\mathcal{A}_{\infty}$-isomorphes (i.e. les structures
$\{\mu_{i}\}_{i\geq 2}$ et $\{\mu'_{i}\}_{i\geq 2}$  sont $\mathcal{A}_{\infty}$-conjugu\'ees).
\end{cor}

Un mot des preuves et des techniques utilis\'ees pour d\'emontrer les th\'eor\`emes
\ref{ti1} et \ref{ti2}. Ces deux r\'esultats utilisent de fa\c{c}on essentielle 
le fait qu'une dg-alg\`ebre propre et lisse $B$ est homotopiquement de pr\'esentation finie
c'est \`a dire que $Map(B,-)$ commute \`a \'equivalence pr\'es avec les colimites filtrantes.  
Ce dernier fait se d\'emontre \`a l'aide de la th\'eorie homotopique des dg-cat\'egories, 
et une preuve se trouve dans \cite{tova1}. 
Le fait que les dg-alg\`ebres propres et lisses soient 
homotopiquement de pr\'esentation finie implique presque imm\'ediatement le point $(2)$ du th\'eor\`eme \ref{ti1}. En effet, 
si on note $(i_{n})_{!}$ l'adjoint \`a gauche du foncteur qui \`a une
dg-alg\`ebre associe sa $\mathcal{A}_{n}$-alg\`ebre sous-jacente, et $\mathbb{L}(i_{n})_{!}$
son foncteur d\'eriv\'e \`a gauche, on montre que l'on a
(voir Prop. \ref{p1})
$$Hocolim_{n} \mathbb{L}(i_{n})_{!}i^{*}_{n}(B) \simeq B.$$
Ainsi, si $B$ est homotopiquement de pr\'esentation finie, il existe un entier $n$ tel que le morphisme naturel
$\mathbb{L}(i_{n})_{!}i^{*}_{n}(B) \longrightarrow B$ poss\`ede une section \`a homotopie pr\`es, 
car l'identit\'e de $B$ se factorise alors par un des $\mathbb{L}(i_{n})_{!}i^{*}_{n}(B)$. 
Ceci implique alors que
$$[B,B'] \longrightarrow [i_{n}^{*}(B),i_{n}^{*}(B')]\simeq [\mathbb{L}(i_{n})_{!}i^{*}_{n}(B),B']$$
poss\`ede une r\'etraction. Nous n'avons pas trouv\'e 
de preuve aussi directe pour le point $(1)$ du th\'eor\`eme \ref{ti1}. La preuve que nous en 
donnons utilise des techniques \'el\'ementaires de \emph{g\'eom\'etrie alg\'ebrique
au-dessus des spectres}, encore appel\'ee \emph{nouvelle g\'eom\'etrie alg\'ebrique
corageuse} (voir \cite{hagII,tova2}). Le point cl\'e est ici de remarquer que pour deux dg-alg\`ebres
propres et lisses $B$ et $B'$ le foncteur $Eq(B,B')$ des \'equivalences entre
$B$ et $B'$, d\'efini sur la cat\'egorie des anneaux en spectres commutatifs, 
en repr\'esentable (i.e. est un sch\'ema affine au-dessus de la cat\'egorie
des spectres, voir Lem. \ref{l9}). Ceci est un fait remarquable de la g\'eom\'etrie alg\'ebrique
au-dessus des spectres, car ce foncteur restreint aux anneaux 
commutatifs usuels n'est en g\'en\'eral pas rep\'esentable par un sch\'ema affine
(sauf si $B$ et $B'$ sont cohomologiquement concentr\'ees en degr\'e $0$). 
De plus, comme $B$ est homotopiquement de pr\'esentation finie 
le foncteur $Eq(B,B')$ est encore homotopiquement de pr\'esentation finie, 
ce qui implique que l'anneau en spectres commutatif qui le repr\'esente
est lui aussi homotopiquement de pr\'esentation finie. On montre alors
de m\^eme que le foncteur $Eq_{n}(B,B')$ des \'equivalences entre les $\mathcal{A}_{n}$-alg\`ebres
$i_{n}^{*}(B)$ et $i_{n}^{*}(B')$ est repr\'esentable par un 
anneau en spectres commutatifs. Enfin, comme nous
l'avons d\'ej\`a vu plus haut on a $Eq(B,B')\simeq Holim_{n}Eq_{n}(B,B')$, 
et le caract\`ere de pr\'esentation finie de $Eq(B,B')$ implique qu'il existe 
un entier $n$ tel que $Eq(B,B') \longrightarrow Eq_{n}(B,B')$ poss\`ede
une r\'etraction. 

Le th\'eor\`eme \ref{ti2} quand \`a lui se d\'eduit du th\'eor\`eme \ref{ti1}
et d'une propri\'et\'e de quasi-compacit\'e du foncteur des classes d'\'equivalences
$k$-dg-alg\`ebres propres, lisses de type $\nu$ et de
dimension cohomologique inf\'erieure \`a $d$. En effet, nous montrons (voir Thm. \ref{t2}) l'existence
d'une $k$-alg\`ebre commutative $A_{0}$ et d'une $A_{0}$-dg-alg\`ebre $B_{0}$
propre, lisse de type $\nu$ et de dimension cohomologique inf\'erieure \`a $d$, telle que
pour toute $k$-dg-alg\`ebre $B$ propre, lisse de type $\nu$ et de dimension cohomologique inf\'erieure \`a $d$
il existe des \'el\'ements $f_{1}, \dots, f_{m} \in k$ avec $\sum f_{i}=1$, et des morphismes
$A_{0} \longrightarrow k[f_{i}^{-1}]$ tels que pour tout $i$ on ait
$$B_{0}\otimes_{A_{0}}^{\mathbb{L}}k[f_{i}^{-1}] \simeq B\otimes_{k}k[f_{i}^{-1}],$$
dans la cat\'egorie homotopique des $k[f_{i}^{-1}]$-dg-alg\`ebres. 
En d'autre termes, le sch\'ema affine $Spec\, A_{0}$ domine, au sens de la topologie
de Zariski, le foncteur des classes d'\'equivalences de dg-alg\`ebres
propres, lisses de type $\nu$ et de dimension cohomologique $d$. Le th\'eor\`eme
\ref{ti1} appliqu\'e aux $A_{0}\otimes_{k}A_{0}$-dg-alg\`ebres $B_{0}\otimes_{A_{0}}^{\mathbb{L}}(A_{0}\otimes_{k}A_{0})$
et $(A_{0}\otimes_{k}A_{0})\otimes_{A_{0}}^{\mathbb{L}}B_{0}$ implique facilement le th\'eor\`eme
\ref{ti2}. \\

Pour terminer cette introduction, signalons que les r\'esultats de ce travail peuvent, et probablement doivent, 
\^etre consid\'er\'es dans le contexte des champs alg\'ebriques sup\'erieurs et m\^eme
d\'eriv\'es (voir par exemple \cite{to}). On peut en effet esp\'erer que le champ
des dg-alg\`ebres propres et lisses est un \emph{$D^{-}$-champ localement g\'eom\'etrique et localement de 
pr\'esentation fini} au sens 
de \cite{tova1}, c'est \`a dire une r\'eunion croissante de n-champs d'Artin (au sens d\'eriv\'e) de type fini. De plus, 
notre caract\'erisation des familles quasi-compactes de dg-alg\`ebres propres et lisses 
laisse penser que le sous-champ des dg-alg\`ebres dont le type et la dimension cohomologique
sont born\'es est un champ d'Artin de type fini. Une approche naturelle pour d\'emontrer le 
caract\`ere alg\'ebrique du champ des dg-alg\`ebres propres et lisses et 
d'appliquer le crit\`ere d'Artin, \'etendu par J. Lurie au cadre d\'eriv\'e. Cependant, la v\'erification
des conditions de ce crit\`ere ne parrait pas totalement \'evidente. Par exemple, montrer que ce champ
est localement de pr\'esentation finie ne semble imm\'ediat (voir le papier compagon \cite{to3} pour une
preuve de ce fait). Il en est de m\^eme de l'alg\`ebrisation des d\'eformations formelles. La question 
du caract\`ere alg\'ebrique du champ des dg-alg\`ebres propres et lisses reste donc
ouverte, et c'est pour cette raison que le point de vue des champs n'a pas \'et\'e adopt\'e dans
ce travail. \\

En effet, les th\'eor\`emes \ref{ti1} et \ref{ti2} sont des 
cons\'equences du fait que le champ des dg-

\bigskip 

\textbf{Remerciements:} Je tiens \`a remercier D. Kaledin et T. Pantev pour des discussions
sur les propri\'et\'es de finitudes des dg-alg\`ebres propres et lisses qui ont 
inspir\'ees les r\'esultats de ce travail. \\

\bigskip

\textbf{Conventions et notations:}
Nous nous excusons de n\'egliger les consid\'erations d'univers, et nous laissons le soins
au lecteur de fixer des univers lorsque cela s'av\`ere n\'ecessaire. 

Notre r\'ef\'erence pour les cat\'egories de mod\`eles est \cite{ho}. Pour une
cat\'egorie de mod\`eles $M$ nous notons $Ho(M)$ sa cat\'egorie homotopique, 
et $[-,-]$ l'ensemble des morphismes dans $Ho(M)$. Pour deux objets
$x$ et $y$ dans $M$ nous notons $Map_{M}(x,y)$, ou  bien
$Map(x,y)$ si le contexte est clair, l'ensemble simplicial des
morphismes tel que d\'efini dans \cite[\S 5]{ho}. L'objet
$Map(x,y)$ sera souvent consid\'er\'e directement dans $Ho(SEns)$.
Nous noterons $Map^{eq}(x,y)$ le sous-ensemble simplicial
de $Map(x,y)$ qui est la r\'eunion des composantes connexes
correspondant aux \'equivalences entre $x$ et $y$. 
Enfin, les produits fibr\'es homotopiques seront not\'e $x\times_{z}^{h}y$. 

Tous les monoides consid\`er\'es dans ce travail seront associatifs et unitaires
(ainsi, tous nos anneaux seront associatifs et unitaires). De m\^eme, tous
les modules sur un monoide seront unitaires (ainsi, tous nos modules sur un anneau 
seront unitaires). Nous mettons en garde le lecteur que pour un anneau $B$
la notation $B-Mod$ fait r\'ef\'erence \`a la cat\'egorie des \emph{$B$-dg-modules}, et 
donc des complexes de $B$-modules, et non pas \`a la cat\'egorie des $B$-modules
au sens usuels (cette derni\`ere ne sera pas not\'ee).

Pour un anneau commutatif $k$ on note $C(k)$ la cat\'egorie des complexes non born\'es
de $k$-modules. On la munit de sa structure de mod\`eles projectives pour laquelle
les \'equivalences sont les quasi-isomorphismes et les fibrations sont 
les \'epimorphismes. 

\section{Dg-alg\`ebres}

Pour tout ce premier chapitre nous fixons un anneau commutatif $k$. \\

Commen\c{c}ons par rappeler que la cat\'egorie des $k$-dg-alg\`ebres  
est par d\'efinition la cat\'egorie des monoides
associatifs et unitaires dans la cat\'egorie monoidale $C(k)$ des complexes
de $k$-modules (non born\'es). Elle sera not\'ee $k-dg-alg$. 

\subsection{La th\'eorie homotopique des dg-alg\`ebres}

On munit la cat\'egorie $k-dg-alg$ de la structure de cat\'egorie de mod\`eles pour laquelle
les fibrations sont les \'epimorphismes et le \'equivalences sont les
quasi-isomorphismes. L'existence de cette structure de mod\`eles se d\'eduit 
des th\'eor\`emes g\'en\'eraux de \cite{ss} appliqu\'es \`a la cat\'egorie
de mod\`eles monoidales $C(k)$ des complexes de $k$-modules. 

On donne ci-dessous quelques propri\'et\'es de la cat\'egorie de mod\`eles des $k$-dg-alg\`ebres. Elles
seront utilis\'ees de fa\c{c}on implicite par la suite. 

\begin{itemize}
\item La cat\'egorie de mod\`eles $k-dg-alg$ est engendr\'ee par cofibration. Pour ensembles g\'en\'erateurs de
cofibrations et de cofibrations triviales on peut prendre l'image de ceux de la cat\'egorie de mod\`eles $C(k)$
par le foncteur ''$k$-dg-alg\`ebre libre'' $C(k) \longrightarrow k-dg-alg$. 

\item Le foncteur d'oubli 
$$k-dg-alg \longrightarrow C(k)$$
pr\'eserve les cofibrations (voir \cite{ss}). En particulier, 
une $k$-dg-alg\`ebre cofibrante est aussi cofibrante comme complexe de $k$-modules, et 
en particulier est un complexe de $k$-modules projectifs sur $k$. 

\item Pour un morphisme d'anneaux commutatifs  $u : k \longrightarrow k'$ 
on dispose d'une adjonction de Quillen
$$-\otimes_{k}k' : k-dg-alg \leftrightarrows k'-dg-alg : f,$$
o\`u l'adjoint \`a droite $f$ est le foncteur d'oubli. L'adjonction d\'eriv\'ee sera not\'ee
$$-\otimes_{k}^\mathbb{L} k' : Ho(k-dg-alg) \leftrightarrows Ho(k'-dg-alg) : f,$$
et en g\'en\'eral nous omettrons de noter le foncteur $f$. 

\item La cat\'egorie de mod\`eles $k-dg-cat$ est \emph{compactement engendr\'ee}
au sens de \cite{tova1}. En particulier, les objets homotopiquement de pr\'esentation finie sont 
exactement les objets \'equivalents aux retractes d'objets cellulaires finis. De fa\c{c}on plus pr\'ecise, 
 pour tout entier $n$ notons $S_{k}(n)$ la $k$-dg-alg\`ebre libre engendr\'ee par un 
\'el\'ement $x$ en degr\'e $-n$ avec $d(x)=0$, et $D_{k}(n)$ la $k$-dg-alg\`ebre libre engendr\'ee
par un \'el\'ement $\alpha$ en degr\'e $-n-1$ et un \'el\'ement $\beta$ en degr\'e $-n$ avec 
$d(\alpha)=\beta$. On dispose d'un morphisme naturel $S_{k}(n) \longrightarrow D_{k}(n+1)$
qui envoie $x$ sur $\beta$. On note $I$ l'ensemble de ces morphismes
$S_{k}(n) \longrightarrow D_{k}(n+1)$ pour $n\in \mathbb{Z}$, qui est un ensemble g\'en\'erateur de cofibrations pour $k-dg-alg$. 

Un objet $B\in k-dg-Akg$ est $I$-cellulaire fini s'il existe une suite finie de diagrammes
cocart\'esiens dans $k-dg-alg$
$$\xymatrix{
B_{i} \ar[r] & B_{i+1} \\
S_{k}(n_{i}) \ar[r] \ar[u] & D_{k}(n_{i}+1), \ar[u]}$$
avec $B_{r}=B$ pour un certain $r$, et $B_{0}=k$. 

Un objet $B$ est homotopiquement de pr\'esentation finie si pour tout
syst\`eme inductif filtrant d'objets $(C_{\alpha})$dans $k-dg-alg$ le morphisme naturel
$$Colim_{\alpha}Map(B,Colim_{\alpha}C_{\alpha}) \longrightarrow Map(B,Colim_{\alpha}C_{\alpha})$$
est une \'equivalence. 

On montre alors que les objets homotopiquement de pr\'esentation finie dans $k-dg-alg$ sont 
exactement les objets \'equivalents \`a des r\'etractes d'objets $I$-cellulaires
finis (voir \cite{tova1}). 

\end{itemize}

\subsection{Cat\'egories d\'eriv\'ees}

Pour une $k$-dg-alg\`ebre $B$, on d\'efinit la cat\'egorie des 
$B$-dg-modules comme \'etant la cat\'egorie des modules \`a gauche dans 
$C(k)$ sur le monoide $B$ (voir \cite{ss}). Cette cat\'egorie sera not\'ee $B-Mod$, et 
est munie de la structure de mod\`eles pour laquelle les fibrations sont les
\'epimorhismes et les \'equivalences sont les quasi-isomorphismes. 

\begin{df}\label{d1}
La \emph{cat\'egorie d\'eriv\'ee} d'une $k$-dg-alg\`ebre $B$ est la cat\'egorie homotopique
de la cat\'egorie de mod\`eles $B-Mod$. Elle est not\'ee
$$D(B):=Ho(B-Mod).$$
\end{df}

La cat\'egorie de mod\`eles $B-Mod$ est stable au sens de \cite[\S 7]{ho}, ainsi la cat\'egorie
$D(B)$ h\'erite d'une structure naturelle de cat\'egorie triangul\'ee (les triangles distingu\'es \'etant les
images dans $D(B)$ des suites exactes de cofibrations dans $B-Mod$).  

Pour un morphisme de $k$-dg-alg\`ebres $B \longrightarrow B'$ on dispose d'une adjonction de Quillen
$$B'\otimes_{B}- : B-Mod \leftrightarrows B'-Mod : f,$$
o\`u l'adjoint \`a droite $f$ est le foncteur d'oubli. Ceci induit une adjonction au niveau des cat\'egories
d\'eriv\'ees
$$B'\otimes^{\mathbb{L}}_{B}-: D(B) \leftrightarrows D(B') : f,$$
et comme pr\'ecemment nous oublierons g\'en\'eralement de noter le foncteur $f$.
Cette adjoint est de plus une \'equivalence de Quillen si le morphisme $B \longrightarrow B'$ est un
quasi-isomorphisme (voir \cite{ss}). 

Rappelons que pour une cat\'egorie triangul\'ee $T$ et un objet $X\in T$ on note
$<X> \subset T$ la sous-cat\'egorie triangul\'ee \'epaisse (i.e. stable par facteurs directs) engendr\'ee par l'objet $X$ 
(voir e.g. \cite{ne}).
Par d\'efinition, $<X>$ est la plus petite sous-cat\'egorie 
triangul\'ee \'epaisse de $T$ contenant $X$.

\begin{df}\label{d2}
La \emph{cat\'egorie d\'eriv\'ee parfaite} d'une $k$-dg-alg\`ebre $B$ est d\'efinie par
$$D_{parf}(B)=<B>\subset D(B).$$
Les objets de $D_{parf}(B)$ seront appel\'es les \emph{$B$-dg-modules parfaits}.
\end{df}

On rappelle que les $B$-dg-modules parfaits poss\`edent les propri\'et\'es suivantes.

\begin{itemize}

\item Un objet de $D(B)$ est parfait si et seulement s'il est compact au sens de \cite{ne}. De m\^eme, un $B$-dg-module
est parfait si et seulement s'il est homotopiquement de pr\'esentation finie dans la cat\'egorie de mod\`eles $B-Mod$.
Ansi, un $B$-dg-module est parfait si et seulement s'il est \'equivalent \`a un r\'etracte d'un $B$-dg-module
$I$-cellulaire fini (voir \cite{tova1} pour plus de d\'etails). 

\item Les $B$-dg-modules sont stables par r\'etractes, extension et d\'ecalage.

\item Lorsque $B$ est une $k$-alg\`ebre (non dg), les $B$-dg-modules parfaits sont 
exactement les complexes de $B$-modules quasi-isomorphes \`a des complexes
born\'es de $B$-modules projectifs de type fini. 

\item Les objets parfaits sont stables par changement de bases: pour un morphisme de $k$-dg-alg\`ebres $B \longrightarrow B'$ 
le foncteur $-\otimes_{B}^{\mathbb{L}}B'$ envoie $B$-dg-modules parfaits sur $B'$-dg-modules
parfaits. 

\end{itemize}

\subsection{Lissit\'e et propret\'e}

Pour deux $k$-dg-alg\`ebres $B$ et $B'$ on peut former leur produit tensoriel
$B\otimes_{k}B'$ qui est encore une $A$-dg-alg\`ebre. Cette construction peut se d\'eriver \`a gauche
en posant
$$B\otimes_{k}^{\mathbb{L}}B':=Q(B)\otimes_{k}B',$$
o\`u $Q$ est un foncteur de remplacement cofibrant dans $k-dg-alg$. Ce produit tensoriel d\'eriv\'e est 
alors compatible avec la notion de quasi-isomorphisme et induit une structure monoidale sym\'etrique
$$-\otimes_{k}^{\mathbb{L}} - : Ho(k-dg-alg) \times Ho(k-dg-alg) \longrightarrow Ho(k-dg-alg).$$
Notons que cette structure monoidale est compatible avec celle sur la cat\'egorie
$D(k)=Ho(C(k))$, dans le sens o\`u le diagramme suivant commute \`a un isomorphisme naturel pr\`es
$$\xymatrix{
Ho(k-dg-alg) \times Ho(k-dg-alg) \ar[r] \ar[d] & Ho(k-dg-alg) \ar[d] \\
D(k) \times D(k) \ar[r] & D(k), }$$
o\`u les foncteurs verticaux sont les foncteurs d'oubli des structures d'alg\`ebres, et les foncteurs
horizontaux sont les produits tensoriels d\'eriv\'es $-\otimes_{k}^{\mathbb{L}}-$.  \\

Pour une $k$-dg-alg\`ebre $B$ on consid\`ere $B$ comme un dg-module sur 
$B\otimes_{k} B^{op}$, o\`u le facteur $B$ op\`ere par multiplication \`a gauche et le facteur
$B^{op}$ op\`ere par multiplication \`a droite. Par ailleurs, le morphisme naturel 
$Q(B) \longrightarrow B$ induit un morphisme de $k$-dg-alg\`ebres
$$B\otimes_{k}^{\mathbb{L}}B^{op} \longrightarrow B\otimes_{k}B^{op},$$
ce qui permet de voir $B$ aussi comme un $B\otimes_{k}^{\mathbb{L}}B^{op}$-dg-module. 
On consid\`ere ainsi $B$ comme un objet dans $D(B\otimes_{k}^{\mathbb{L}}B^{op})$. 

\begin{df}\label{d3}
Soit $B$ une $k$-dg-alg\`ebre. 
\begin{enumerate}
\item La $k$-dg-alg\`ebre $B$ est \emph{propre} si 
$B$ est parfait comme objet de $D(k)$ (i.e. si le complexe de $k$-modules
sous-jacent \`a $B$ est parfait). 

\item La $k$-dg-alg\`ebre $B$ est \emph{lisse} si 
$B$ est parfait comme objet dans $D(B\otimes_{k}^{\mathbb{L}}B^{op})$.

\end{enumerate}
\end{df}

\begin{rmk}\label{r1}
\begin{enumerate}
\item \emph{Il est important de remarquer que les notions pr\'ec\'edentes
de lissit\'e et de propret\'e d\'ependent du choix de l'anneaux de base $k$.}

\item \emph{Pour plus d'explication sur la terminologie de} lisse \emph{et} propre \emph{on 
renvoie \`a \cite{koso,tova1}, o\`u le lecteur trouvera aussi des exemples de telles dg-alg\`ebres}.
\end{enumerate}
\end{rmk}

Un r\'esultat cl\'e concernant les dg-alg\`ebres propres et lisses est le th\'eor\`eme suivant, tir\'e
de \cite{tova1}, et dont la preuve, qui utilise la th\'eorie homotopique
des dg-cat\'egories, ne sera pas reproduite ici.  Ce th\'eor\`eme est \`a la base du th\'eor\`eme de finitude que nous
d\'emontrerons au pargraphe 4.

\begin{thm}{(\cite{tova1})}\label{t1}
Toute $k$-dg-alg\`ebre propre et lisse est homotopiquement de pr\'esentation finie. 
\end{thm}

\subsection{Dimension cohomologique des dg-alg\`ebres lisses}

Soit $B$ une $k$-dg-alg\`ebre cofibrante. On consid\`ere
l'adjonction de Quillen
$$B\otimes_{k} - : C(k) \leftrightarrows B-Mod : f,$$
o\`u $f$ est le foncteur d'oubli. Comme $B$ est cofibrante elle est plate, et  ces deux foncteurs 
pr\'eservent les \'equivalences. 
En posant $\mathcal{F}:=(B\otimes_{k} -)\circ f$, cette adjonction
d\'efinit un foncteur
$$\mathcal{R}_{*} : B-Mod \longrightarrow sB-Mod,$$
o\`u $\mathcal{R}_{*}(M)$ est l'objet simplicial d\'efini par
$$\begin{array}{cccc}
\mathcal{R}_{*}(M) : & \Delta^{op} & \longrightarrow & B-Mod \\
 & [n] & \mapsto & \mathcal{F}^{\circ n}(M).
\end{array}$$
Pour tout $M\in B-Mod$, l'objet simplicial $\mathcal{R}_{*}(M)$ est augment\'e sur $M$, 
et n'est autre que la r\'esolution simpliciale standard de $M$ associ\'ee \`a l'adjonction pr\'ec\'edente
(voir \cite{il}). En particulier, le morphisme induit
$$|\mathcal{R}_{*}(M)|:= Hocolim_{[n]\in \Delta^{op}}\mathcal{R}_{n}(M) \longrightarrow M$$
est un isomorphisme dans $D(B-Mod)$ pour tout $B$-dg-module $M$. 
De plus, pour $M\in B-Mod$, l'objet $|\mathcal{R}_{*}(M)|$ est naturellement \'equivalent \`a la colimite filtrante
$$|\mathcal{R}_{*}(M)|\simeq Colim_{k\in\mathbb{N}} |Sq_{k}\mathcal{R}_{*}(M)|,$$
o\`u $|Sq_{k}\mathcal{R}_{*}(M)|$ est la colimite homotopique du diagramme
$\mathcal{R}_{*}(M)$ restreint \`a la sous-cat\'egorie pleine $\Delta_{\leq k}^{op} \subset \Delta^{op}$ form\'ee des 
objets $[n]$ avec  $n\leq k$.  On obtient ainsi un isomorphisme dans $D(B)$
$$Colim_{k\in\mathbb{N}} |Sq_{k}\mathcal{R}_{*}(M)| \simeq M.$$

Soit maintenant $B$ une $k$-dg-alg\`ebre quelconque, et $Q(B) \longrightarrow B$
un mod\`ele cofibrant. Pour un $B$-dg-module $M$, on peut consid\'erer
$M$ comme un $Q(B)$-dg-module et appliquer la construction pr\'ec\'edente.
Cela nous fournit alors un isomorphisme dans $D(Q(B))$
$$Colim_{k\in\mathbb{N}} |Sq_{k}\mathcal{R}_{*}(M)| \simeq M,$$
o\`u la r\'esolution $\mathcal{R}_{*}(M)$ est calcul\'ee dans les $Q(B)$-dg-modules. 
Comme $-\otimes_{Q(B)}^{\mathbb{L}}B : D(Q(B)) \longrightarrow D(B)$ est une
\'equivalence de cat\'egorie, nous consid\`ererons aussi cet isomorphisme
dans $D(B)$ (sans pour
autant sp\'ecifier le changement de bases $-\otimes_{Q(B)}^{\mathbb{L}}B$). Nous mettons
en garde le lecteur que l'objet $|\mathcal{R}_{*}(M)|$ dans $D(Q(B))$, vu
comme un objet simplicial dans $D(B)$ \`a l'aide du foncteur $-\otimes_{Q(B)}^{\mathbb{L}}B$
n'est pas la r\'esolution libre standard du $B$-dg-module $M$ en g\'en\'eral.

\begin{df}\label{d4}
\begin{enumerate}
\item 
Un dg-module sur une $k$-dg-alg\`ebre $B$ est \emph{de dimension cohomologique inf\'erieure \`a $d\in \mathbb{N}$} 
si le morphisme d'augmentation
$$|Sq_{d}\mathcal{R}_{*}(M)| \longrightarrow M$$
poss\`ede une section dans $D(Q(B))\simeq D(B)$.
\item Un dg-module sur une $k$-dg-alg\`ebre est \emph{de dimension cohomologique infinie} s'il n'est 
pas de dimension cohomologique inf\'erieure \`a $d$ pour tout $d\in \mathbb{N}$. Il est
\emph{de dimension cohomologique finie sinon}.
\item Une $k$-dg-alg\`ebre $B$ est \emph{de dimension cohomologique inf\'erieure \`a $d\in \mathbb{N}$}
si le $B\otimes_{k}^{\mathbb{L}}B^{op}$-dg-module $B$ est de dimension cohomologique inf\'erieure
\`a $d$.
\item Une  $k$-dg-alg\`ebre est \emph{de dimension cohomologique infinie} si elle n'est 
pas de dimension cohomologique inf\'erieure \`a $d$ pour tout $d\in \mathbb{N}$. Elle est
\emph{de dimension cohomologique finie sinon}.
\end{enumerate}
\end{df}

\begin{rmk}\label{r2}
\begin{enumerate}
\item
\emph{Commen\c{c}ons pas remarquer que la notion de dimension cohomologique d'une
$k$-dg-alg\`ebre $B$ d\'epend du choix de la base $k$ (car 
$B\otimes_{k}^{\mathbb{L}}B^{op}$ d\'epend de cette base).}
\item \emph{Il est facile de voir qu'un dg-module qui est de dimension cohomologique inf\'erieure
\`a $d$ est aussi de dimension cohomologique inf\'erieure \`a $d'$ pour tout $d'\geq d$. 
De m\^eme pour une $k$-dg-alg\`ebre.}
\item \emph{Remarquons aussi que pour tout  $B$-dg-module $M$ on a
$$M\simeq Colim_{k\in \mathbb{N}} |Sq_{k}\mathcal{R}_{*}(M)|.$$
Ainsi, si $M$ est parfait il existe un entier $k\geq 0$ tel que le morphisme d'augmentation
$$|Sq_{k}\mathcal{R}_{*}(M)| \longrightarrow M$$
poss\`ede une section dans $D(B)$. Ceci montre qu'un tout
dg-module parfait est de dimension cohomologique finie. En particulier, 
une $k$-dg-alg\`ebre lisse est toujours de dimension cohomologique finie. 
R\'eciproquement on peut voir qu'un $B$-dg-module $M$ qui est parfait sur 
$k$ et qui est de dimension cohomologique finie est aussi un $B$-dg-module parfait. Ainsi, 
une $k$-dg-alg\`ebre propre est lisse si et seulement si elle est 
de dimension cohomologique finie.}
\end{enumerate}
\end{rmk}

\subsection{Dg-alg\`ebres, $\mathcal{A}_{\infty}$-Alg\`ebres et $\mathcal{A}_{n}$-alg\`ebres}

Commen\c{c}ons par rappeler la notion d'op\'erade dans une
cat\'egorie monoidale sym\'etrique (voir \cite{sp}). Dans ce qui suit
nous ne consid\`ererons que des op\'erades ''non-$\Sigma$" et 
unitaires, et les alg\`ebres sur ces op\'erades seront toujours unitaires.
D'apr\`es \cite{sp} la cat\'egorie des op\'erades $Op(C(k))$ dans 
la cat\'egorie des complexes de $k$-modules est munie d'une structure
de cat\'egorie de mod\`eles pour laquelle les fibrations et les \'equivalences
sont d\'efinis sur les complexes sous-jacents. 

On rappelle d'apr\`es l'existence des op\'erades $\mathcal{A}_{\infty}$ et 
$\mathcal{A}_{n}$ au-dessus dans la cat\'egorie monoidale sym\'etrique $C(k)$. 
Rappelons, qu'en tant  qu'op\'erade dans les $k$-modules gradu\'es (i.e. en 
oubliant la diff\'erentielle) $\mathcal{A}_{\infty}$ est librement engendr\'ee
par des op\'erations $\mu_{i}\in \mathcal{A}_{\infty}(i)$ de degr\'es $2-i$. 
En utilisant les notations de \cite{ma} on a 
$$\mathcal{A}_{\infty}=\Gamma(\mu_{2},\mu_{3}, \dots, \mu_{i}, \dots).$$
La diff\'erentielle sur $\mathcal{A}_{\infty}$ est alors d\'efinie par la 
formule suivante
$$\partial (\mu_{n})=\sum_{i+j=n+1, i,j\geq 2} \sum_{s=0}^{s=n-1}(-1)^{j+s(j+1)}\mu_{i}
(\mathbf{1}^{\otimes s}\otimes \mu_{j} \otimes \mathbf{1}^{\otimes i-s-1}).
$$

L'op\'erade $\mathcal{A}_{\infty}$ est cofibrante pour la structure de mod\`eles
de \cite{sp}, et les alg\`ebres sur $\mathcal{A}_{\infty}$ sont pr\'ecis\`emment les
$\mathcal{A}_{\infty}$-alg\`ebres au sens usuel (e.g. de \cite{ke1,lef}). Dapr\`es \cite{sp}, la cat\'egorie
$k-\mathcal{A}_{\infty}-alg$, des alg\`ebres unitaires sur $\mathcal{A}_{\infty}$ est munie
d'une structure de cat\'egorie de mod\`eles pour laquelle les fibrations et les 
\'equivalences sont d\'efinis sur les complexes sous-jacents. De plus, 
il existe un morphisme naturel d'op\'erades $p : \mathcal{A}_{\infty} \longrightarrow Ass$, 
o\`u $Ass$ est l'op\'erade des alg\`ebres associatives. Ce morphisme induit une
adjonction de Quillen
$$p_{!} : k-\mathcal{A}_{\infty}-alg \leftrightarrows k-dg-alg : p^{*}$$
qui s'av\`ere \^etre une \'equivalence de Quillen. A travers ectte \'equivalence de Quillen nous
identifierons souvent la th\'eorie homotopique des $k$-dg-alg\`ebres
avec  celle de $k$-$\mathcal{A}_{\infty}$-alg\`ebres, sans pour autant mentionner
l'\'equivalence $(\mathbb{L}p _{!},p^{*})$. 

Pour un entier $n\geq 2$, on note $\mathcal{A}_{n} \subset \mathcal{A}_{\infty}$ la sous-op\'erade
engendr\'ee par les op\'erations $\mu_{i}$ pour $2\leq i\leq n$. Les op\'erades 
$\mathcal{A}_{n}$ sont cofibrantes, et les alg\`ebres au-dessus de $\mathcal{A}_{n}$
sont pr\'ecis\`emment les $\mathcal{A}_{n}$-alg\`ebres. Les morphismes
d'inclusions $i_{n} : \mathcal{A}_{n} \hookrightarrow \mathcal{A}_{\infty}$
induisent des adjonctions de Quillen
$$(i_{n})_{!} : k-\mathcal{A}_{n}-alg \leftrightarrows k-\mathcal{A}_{\infty}-alg :i_{n}^{*}.$$

\begin{lem}\label{l1}
Le morphisme naturel 
$$Hocolim_{n} \mathcal{A}_{n} \longrightarrow \mathcal{A}_{\infty}$$
est une \'equivalence dans $Op(C(k))$. 
\end{lem}

\textit{Preuve:} Cela se d\'eduit du fait que les colimites filtrantes sont 
aussi des colimites homotopiques dans $Op(C(k))$, et du fait 
\'evident que $Colim_{n}\mathcal{A}_{n}\simeq \mathcal{A}_{\infty}$.
\hfill $\Box$ \\

Nous allons d\'eduire du lemme \ref{l1} le r\'esultat suivant.

\begin{prop}\label{p1}
Soient $B$ et $B'$ deux $k$-dg-alg\`ebres. Alors le morphisme naturel d'ensembles simpliciaux
$$Map_{k-dg-alg}(B,B') \longrightarrow Map_{k-\mathcal{A}_{\infty}}(p^{*}(B),p^{*}(B')) \longrightarrow
Holim_{n}Map_{k-\mathcal{A}_{n}-alg}(i_{n}^{*}p^{*}(B),i_{n}^{*}p^{*}(B'))$$
est un \'equivalence. 
\end{prop}

\textit{Preuve:} Le premier de ces morphismes est une \'equivalence car 
$(p_{!},p^{*})$ est une \'equivalence de Quillen. Il nous faut donc montrer que
le second morphisme est aussi une \'equivalence. 

Pour une cat\'egorie de mod\`eles $M$ nous noterons $|M|$
l'ensemble simplicial nerf de la sous-cat\'egorie des \'equivalences dans $M$. 
On consid\`ere alors les morphismes naturels
$$\xymatrix{
 |k-\mathcal{A}_{\infty}-alg| \ar[d]^-{q} \ar[r] & Holim_{n}|k-\mathcal{A}_{n}-alg| \ar[dl]^-{r} \\
  |C(k)|, & }$$
o\`u les morphismes vers $|C(k)|$ sont induits par les foncteurs qui oublient les structures
d'alg\`ebres et ne retiennent que les complexes sous-jacents. On utilise alors
les r\'esultats de \cite{re}, ou plut\^ot ses g\'en\'eralisations imm\'ediates au-dessus de $C(k)$.
D'apr\`es \cite[Thm. 1.2.15]{re}, pour $E\in |C(k)|$ un complexe, 
les morphismes induits sur les fibres homotopiques des morphismes $q$ et $r$
sont \'equivalents aux morphismes naturels
$$\xymatrix{
Map_{Op(C(k))}(\mathcal{A}_{\infty},\mathbb{R}\underline{End}(E)) \longrightarrow 
Holim_{n}Map_{Op(C(k))}(\mathcal{A}_{n},\mathbb{R}\underline{End}(E)),}$$
o\`u $\mathbb{R}\underline{End}(E)$ est l'op\'erade d'endomorphismes
d'un mod\`eles cofibrant pour $E$. D'apr\`es le lemme \ref{l1}, ce dernier morphisme
est un \'equivalence. On en d\'eduit donc que le morphisme
$$|k-\mathcal{A}_{\infty}-alg|  \longrightarrow Holim_{n}|k-\mathcal{A}_{n}-alg|$$
est une \'equivalence. 

Notons maintenant $Mor(C(k))$ la cat\'egorie de mod\`eles des morphismes
dans $C(k)$ (munie par exemple de sa structure projective pour laquelle les
fibrations et les \'equivalences sont d\'efinies sur les objets sous-jacents dans $C(k)$). 
La cat\'egorie de mod\`eles $Mor(C(k))$ reste une cat\'egorie de mod\`eles
monoidale sym\'etrique pour laquelle la structure monoidale est d\'efinie objets par 
objets dans $M$. On dipose donc de cat\'egorie de mod\`eles
de $\mathcal{A}_{\infty}$ et $\mathcal{A}_{n}$-alg\`ebres dans $Mor(C(k))$, 
qui s'identifie aux cat\'egories de mod\`eles $Mor(k-\mathcal{A}_{\infty}-alg)$
et $Mor(k-\mathcal{A}_{n}-alg)$.
On peut donc de nouveau appliquer \cite[Thm. 1.2.15]{re} au-dessus de la cat\'egorie
de base $Mor(C(k))$, et par le m\^eme argument que ci-dessus on obtenir une \'equivalence de nerfs
$$|Mor(k-\mathcal{A}_{\infty}-alg)| \longrightarrow Holim_{n}|Mor(k-\mathcal{A}_{n}-alg)|.$$
Pour terminer la preuve de la proposition \ref{p1}, soient $B$ et $B'$ deux points
dans $|k-\mathcal{A}_{\infty}-alg|$. On dispose d'un diagramme homotopiquement 
commutatif
$$\xymatrix{
|Mor(k-\mathcal{A}_{\infty}-alg)|  \ar[r]^-{u} \ar[d] & |k-\mathcal{A}_{\infty}-alg|\times |k-\mathcal{A}_{\infty}-alg|     \ar[d] \\
Holim_{n} |Mor(k-\mathcal{A}_{n}-alg)| \ar[r]_-{v} & Holim_{n} (|k-\mathcal{A}_{n}-alg|\times |k-\mathcal{A}_{n}-alg|), }$$
o\`u les morphismes horizontaux associent \`a un morphisme sa source et son but. 
Les morphismes verticaux \'etant des \'equivalences on en d\'eduit une \'equivalence 
entre les fibres homotopiques de $u$ et de $v$ prises en le point $(B,B')$. Mais d'apr\`es \cite[Thm.8.3]{re2} ce morphisme
induit sur les fibres homotopiques est \'equivalent au morphisme
$$Map_{k-\mathcal{A}_{\infty}-alg}(B,B') \longrightarrow Holim_{n}Map_{k-\mathcal{A}_{n}-alg}(i _{n}^{*}(B),i_{B}^{*}(B')).$$
Ceci termine la preuve de la proposition \ref{p1}. \hfill $\Box$ \\

\begin{cor}\label{c1}
Pour toute $k$-dg-alg\`ebre $B$ (consid\'er\'ee aussi comme une $k-\mathcal{A}_{\infty}$-alg\`ebre) le morphisme naturel
$$Colim_{n} \mathbb{L}(i_{n})_{!}i_{n}^{*}(B) \longrightarrow B$$
est une \'equivalence. 
\end{cor}

\textit{Preuve:}  D'apr\`es la proposition \ref{p1} pour tout $B'\in k-dg-alg$ le morphisme naturel
$$Map_{k-dg-alg}(B,B') \longrightarrow 
Map_{k-dg-alg}(Colim_{n} \mathbb{L}(i_{n})_{!}i_{n}^{*}(B),B') $$
$$\simeq
Holim_{n}Map_{k-dg-alg}(\mathbb{L}(i_{n})_{!}i^{*}_{n}(B),B')
\simeq Holim_{n} Map_{k-\mathcal{A}_{n}-alg}(i^{*}_{n}(B),i_{n}^{*}(B'))$$
est une \'equivalence. Le lemme de Yoneda pour la cat\'egorie
homotopique $Ho(k-dg-alg)$ implique donc que le morphisme
$Colim_{n} \mathbb{L}(i_{n})_{!}i_{n}^{*}(B) \longrightarrow B$ est un 
isomorphisme dans $Ho(k-dg-alg)$ et donc une \'equivalence. 
\hfill $\Box$ \\

\begin{cor}\label{c2}
Pour toute $k$-dg-alg\`ebre $B$ qui est homotopiquement de pr\'esentation finie (e.g. propre et lisse), il existe un 
entier $n$ tel que le morphisme
$$\mathbb{L}(i_{n})_{!}i_{n}^{*}(B) \longrightarrow B$$
poss\`ede une section dans $Ho(k-dg-alg)$. 
\end{cor}

\textit{Preuve:} C'est une cons\'equence imm\'ediate du corollaire \ref{c1} et 
de la d\'efinition d'\^etre homotopiquement de pr\'esentation finie. \hfill $\Box$ \\

\section{Familles quasi-compactes de dg-alg\`ebres propres et lisses}

On fixe un anneau commutatif de base $k$.

\subsection{Faisceaux quasi-compacts}

Notons $\widehat{k-Aff}$ la cat\'egorie des foncteurs $k-Aff^{op}=k-CAlg \longrightarrow Ens$, 
de la cat\'egorie des $k$-alg\`ebres commutatives vers celle des ensembles. Nous
aurons aussi \`a consid\'erer la cat\'egorie $SPr(k-Aff)$, des objets simpliciaux
dans $\widehat{k-Aff}$ (c'est aussi la cat\'egorie des foncteurs $k-CAlg \longrightarrow SEns$).
La cat\'egorie $SPr(k-Aff)$ est munie de sa structure de mod\`eles projective niveaux
par niveaux (i.e. fibrations et \'equivalences sont d\'efinies objets par objets au-dessus
de $k-Aff$, nous n'inlcuons pas la topologie de Zariski dans cette
structure de mod\`eles).  On dispose d'un foncteur
$$\pi_{0} : SPr(k-Aff) \longrightarrow \widehat{k-Aff}$$
qui \`a $F$ associe le pr\'efaisceau $A \mapsto \pi_{0}(F(A))$. Le foncteur
$\pi_{0}$ est adjoint \`a gauche du foncteur d'inclusion 
$\widehat{k-Aff} \longrightarrow SPr(k-Aff)$ qui voit un pr\'efaisceau d'ensembles
comme un pr\'efaisceau d'ensembles simpliciaux constant dans la direction 
simpliciale. 

\begin{df}\label{d5}
Nous dirons qu'un foncteur $F\in  \widehat{k-Aff}$ est \emph{quasi-compact} 
s'il existe un sch\'ema affine
$X$ et un morphisme de pr\'efaisceaux $X \longrightarrow F$ qui induise 
un \'epimorphisme sur les faisceaux associ\'es pour la topologie de Zariski sur $k-Aff$.
\end{df}

\begin{rmk}
\emph{Un morphisme de pr\'efaisceaux $F \longrightarrow G$ qui induit un 
\'epimorphisme sur les faisceaux Zariski associ\'es sera appel\'e un}
\'epimorphisme Zariski local. 
\end{rmk}

En d'autres termes, le pr\'efaisceau $F$ est quasi-compact si et seulement s'il existe
une $k$-alg\`ebre commutative $A_{0}$ et un \'el\'ement $x\in F(A_{0})$, tels que 
pour tout $A\in k-CAlg$ et tout $y\in F(A)$, il existe des \'el\'ements
$a_{1}, \dots, a_{n} \in A$ avec $\sum a_{i}=1$, et des morphismes
$u_{i} : A_{0} \longrightarrow A[a_{i}^{-1}]$, avec $u_{i}(x)=f_{i}(y)$ dans $F(A[a_{i}]^{-1})$, o\`u 
$f_{i} : A \longrightarrow A[a_{i}^{-1}]$ est le morphisme canonique. \\

Le lemme suivant fournit deux proc\'ed\'es pour v\'erifier qu'un 
pr\'efaisceau est quasi-compact.

\begin{lem}\label{l2}
\begin{enumerate}
\item Soit $F \longrightarrow G$ un morphisme dans $\widehat{k-Aff}$. On suppose que
$G$ est quasi-compact et que pour tout sch\'ema affine $X$ et tout morphisme
$X \longrightarrow G$ le pr\'efaisceau $F\times_{G}X$ est quasi-compact. Alors
$F$ est quasi-compact.
\item Soit $F$ et $G$ deux pr\'efaisceaux simpliciaux sur $k-Aff$
$$F,G : k-Aff^{op} \longrightarrow SEns,$$
et $f : F \longrightarrow G$ un morphisme. On suppose que le pr\'efaisceau
$\pi_{0}(F)$ est quasi-compact, et que pour tout sch\'ema affine $X$ et tout morphisme
$X \longrightarrow G$ le pr\'efaisceau $\pi_{0}(F\times_{G}^{h}X)$
est quasi-compact. Alors le pr\'efaisceau 
$\pi_{0}(F)$ est quasi-compact. 
\item Soit $F$ un pr\'efaisceau tel qu'il existe un sch\'ema quasi-compact $X$ 
et un \'epimorphisme Zariski local $X \longrightarrow F$, alors $F$ est quasi-compact.
\end{enumerate}
\end{lem}

\textit{Preuve:} $(1)$ On choisit un \'epimorphisme Zariski local $X \longrightarrow G$, avec $X$ un sch\'ema
affine. Par hypoth\`ese on peut choisir un \'epimorphisme Zariski local $Y \longrightarrow F\times_{G}X$
avec $Y$ un sch\'ema affine.
Le morphisme compos\'e $Y \longrightarrow F\times_{G}X \longrightarrow F$
est un encore un \'epimorphisme Zariski local, car compos\'e de deux \'epimorphismes
Zariski locaux (on utilise ici que les \'epimorphismes Zariski locaux sont 
stables par compositions et changement de bases). \\

$(2)$ Pour tout sch\'ema affine $X$ et tout morphisme $X \longrightarrow G$, 
le morphisme naturel 
$$\pi_{0}(F\times_{G}^{h}X) \longrightarrow \pi_{0}(F)\times_{\pi_{0}(G)}X$$
est un \'epimorphisme de pr\'efaisceaux, et en particulier un 
\'epimorphisme Zariski local. Par hypoth\`ese $\pi_{0}(F\times_{G}^{h}X)$ est 
quasi-compact, et ainsi le pr\'efaisceau $\pi_{0}(F)\times_{\pi_{0}(G)}X$
est aussi quasi-compact. On peut alors appliquer le point $(1)$
au morphisme $\pi_{0}(F) \longrightarrow \pi_{0}(G)$. \\

$(3)$ Le sch\'ema est recouvert par un nombre fini d'ouverts Zariski affines $U_{i}$. On pose
$U:=\coprod_{i}U$ (coproduit pris dans les sch\'emas), 
qui est un sch\'ema affine. Alors, le morphisme compos\'e
$U \longrightarrow X \longrightarrow F$ est un \'epimorphisme Zariski local.  
\hfill $\Box$ \\

\begin{df}\label{d6}
Soit $G \in \widehat{Aff}$ un pr\'efaisceau.
\begin{enumerate}
\item Une \emph{famille d'objets dans $G$} est un sous-ensemble 
$\mathcal{F} \subset G(k)$. 
\item Une famille $\mathcal{F}$ d'objets de $G$ est \emph{quasi-compacte} s'il existe
un sous-pr\'efaisceau quasi-compact $F \subset G$ tel que $\mathcal{F}\subset F(k)$. 
\end{enumerate}
\end{df}

\begin{rmk}
\emph{Il est important de noter que la propri\'et\'e de quasi-compacit\'e d\'epend
d'un pr\'efaisceau ambient $G$. En g\'en\'eral si $G\subset G'$ est une inlcusion
de pr\'efaisceaux, une famille $\mathcal{F}$ d'objets dans $G$ peut tout \`a fait \^etre 
quasi-compacte comme famille d'objets dans $G'$ mais pas comme famille
d'objets dans $G$. Par exemple, on peut prendre pour $G$ un sch\'ema affine
qui n'est pas de type fini sur $k$, pour $G'$ un ouvert non quasi-compact de $G$ et 
$\mathcal{F}:=G'(k)$. La famille $\mathcal{F}$ est quasi-compacte dans $G'$ car 
$G'$ est un sch\'ema affine, mais n'est en g\'en\'eral pas quasi-compacte dans $G$.}
\end{rmk}

\subsection{Caract\'erisation des familles quasi-compactes de dg-alg\`ebres propres et lisses}

Pour tout $A\in k-CAlg$, nous noterons 
$dga(A)$ l'ensemble des classes d'isomorphismes dans la cat\'egorie
$Ho(A-dg-alg)$. Pour un morphisme $A \longrightarrow A'$ dans $k-CAlg$, le foncteur
de changement de bases $A'\otimes_{A}^{\mathbb{L}} - : Ho(A-dg-alg) \longrightarrow
Ho(A'-dg-alg)$ induit une application 
$$dga(A) \longrightarrow dga(A').$$
Ceci d\'efinit un pr\'efaisceau $dga \in \widehat{k-Aff}$. Nous d\'efinissons aussi
$dg-alg^{p,l} \subset dg-alg$ comme \'etant le sous-pr\'efaisceau des
dg-alg\`ebres propres et lisses. 

\begin{df}\label{d6'}
\begin{enumerate}
\item Une \emph{famille de $k$-dg-alg\`ebres propres et lisses} est 
une famille d'objets dans le pr\'efaisceau $dg-alg^{p,l}$ (i.e. un
sous-ensemble de l'ensemble des classes de quasi-isomorphismes
de $k$-dg-alg\`ebres propres et lisses).
\item Une famille $\mathcal{F}$ de $k$-dg-alg\`ebres propres et lisses est \emph{quasi-compacte} si
elle l'est comme famille d'objets de $dg-alg^{p,l}$ (au sens de la d\'efinition \ref{d6}).
\end{enumerate}
\end{df}

En termes plus explicites: une famille de $k$-dg-alg\`ebres propres et lisses $\mathcal{F}$ est quasi-compacte s'il existe
$A\in k-CAlg$ et $B_{0}$ une $A$-dg-alg\`ebre \emph{propre et lisse sur $A$}, 
tels que pour tout $B\in \mathcal{F}$, il existe des \'el\'ements 
$f_{1}, \dots, f_{n} \in k$, avec $\sum f_{i}=1$, et des morphismes
$A \longrightarrow k[f_{i}^{-1}]$, tels que 
$$B\otimes^{\mathbb{L}}_{k}k[f_{i}^{-1}] \simeq B_{0}\otimes^{\mathbb{L}}_{k}k[f_{i}^{-1}]$$
pour tout $i$ (il s'agit d'isomorphismes dans $Ho(k[f_{i}^{-1}]-dg-alg)$). On voit en particulier que lorsque 
$k$ est un anneau local, on peut prendre $n=1$ et $f_{1}=1$. On dispose dans ce cas
d'un $A$-dg-alg\`ebre $B_{0}$ propre et lisse sur $A$, tel que tout $B\in \mathcal{F}$ soit de la forme
$B_{0}\otimes^{\mathbb{L}}_{A}k$ pour un certain morphisme $A\longrightarrow k$. \\

Pour \'enoncer notre th\'eor\`eme de caract\'erisation des
familles quasi-compactes de $k$-dg-alg\`ebres propres et lisse 
nous introduisons la notion de type d'une
$k$-dg-alg\`ebre propre. 

\begin{df}\label{d7}
\begin{enumerate}
\item Un \emph{type} est une application $\nu : \mathbb{Z} \longrightarrow \mathbb{N}$
qui est nulle en dehors d'un intervalle fini (i.e. \`a support fini). 
\item Soit $\nu$ un type et $A$ un anneau commutatif. 
Nous dirons qu'un complexe parfait $E \in D_{parf}(A)$  \emph{est de type
$\nu$ (relativement \`a $A$)} si pour tout corps $K$ et tout morphisme d'anneaux $A \longrightarrow K$ on a 
$$Dim_{K}H^{i}(E\otimes_{k}^{\mathbb{L}}K) \leq \nu(i) \qquad \forall \; i \in \mathbb{Z}.$$
\item Soit $\nu$ un type et $A$ un anneau commutatif. 
Nous dirons qu'un $A$-dg-alg\`ebre propre $B$ \emph{est de type $\nu$}
si son complexe sous-jacent est de type $\nu$. 
\end{enumerate}
\end{df}

Le th\'eor\`eme suivant caract\'erise les familles quasi-compactes de dg-alg\`ebres propres et lisses
comme les familles dont les types et les dimensions cohomologiques sont born\'es. 

\begin{thm}\label{t2}
Soit $\mathcal{F}\subset dg-alg^{p,l}(k)$ une famille de $k$-dg-alg\`ebres propres et lisses.
Alors, $\mathcal{F}$ est quasi-compacte si et seulement s'il existe un type $\nu$ et un entier $d$ tels que les
deux conditions suivante soient satisfaites. 
\begin{enumerate}
\item Toute dg-alg\`ebre $B\in \mathcal{F}$ est de type $\nu$ (relativement \`a $k$). 
\item Pour toute dg-alg\`ebre $B \in \mathcal{F}$, il existe $f_{1}, \dots, f_{n}\in k$ avec 
$\sum f_{i}=1$, tels que chaque $B\otimes^{\mathbb{L}}_{k}k[f_{i}^{-1}]$ soit de dimension cohomologique 
inf\'erieure \`a $d$ (relativement \`a $k[f_{i}^{-1}]$).
\end{enumerate}
\end{thm}

La d\'emonstation de ce th\'eor\`eme va prendre un certain temps, et nous y consacrons
la fin de cette section. La n\'ecessit\'e de la condition est facile et est laiss\'ee au lecteur. 
Nous nous concentrerons sur la suffisance. \\

Soit $\nu : \mathbb{Z} \longrightarrow \mathbb{N}$ un type. Pour toute $k$-alg\`ebre
commutative $A\in k-CAlg$, notons $\mathcal{E}_{\nu}(A)$
le sous-ensemble de l'ensemble des classes d'isomorphismes
de $D_{parf}(A)$, form\'e des complexes parfaits et de type
$\nu$ (sur $A$). Pour un morphisme $A \longrightarrow B$ dans $k-CAlg$, on dispose d'un 
changement de base 
$$B\otimes_{A}^{\mathbb{L}} - : D_{parf}(A) \longrightarrow D_{parf}(B).$$
Ce foncteur pr\'eserve les objets
de type $\nu$. On obtient donc une application
$$B\otimes_{A}^{\mathbb{L}} - : \mathcal{E}_{\nu}(A) \longrightarrow \mathcal{E}_{\nu}(B), $$
qui permet de voir $A \mapsto \mathcal{E}_{\nu}(A)$ comme un objet de
$\widehat{k-Aff}$. 

\begin{lem}\label{l3}
Le pr\'efaisceau $\mathcal{E}_{\nu}$ est quasi-compact.
\end{lem}

\textit{Preuve:} Notons $\mathcal{E}$ le pr\'efaisceau sur $k-Aff$ qui \`a $A\in k-CAlg$ associe 
l'ensemble des classes d'isomorphismes d'objets dans $D_{parf}(A)$. 
Notons
$X \in \widehat{k-Aff}$ qui \`a $A$ associe 
l'ensemble des complexes de $A$-modules
$$\xymatrix{
\dots \ar[r] A^{n_{i+1}} \ar[r]^-{d_{i+1}} & A^{n_{i}} \ar[r]^-{d_{i}} & A^{n_{i-1}} \ar[r]^-{d_{i-1}} & \dots},$$
avec $0\leq n_{i}\leq \nu(i)$ pour tout $i$. Il est facile de voir que  $X$ est repr\'esent\'e par 
une r\'eunion disjointe finie de sch\'emas affines, dont les composantes
param\'etrisent les complexes pour lesquels les $n_{i}$ sont tous fix\'es. 
On dispose d'un morphisme
de pr\'efaisceaux
$$X \longrightarrow \mathcal{E}$$
qui \`a un complexes de $A$-modules comme ci-dessus associe 
sa classe d'isomorphisme dans $D_{parf}(A)$. Il est facile de voir que ce morphisme
se factorise \`a travers le faisceau associ\'e 
$$\xymatrix{
X \ar[r] \ar[d] & a(X) \ar[dl] \\
\mathcal{E} & }$$
(on remarquera que le foncteur $\mathcal{E} : k-CAlg \longrightarrow Ens$ commute
avec les produits finis, et donc  $Hom(\coprod X_{i},\mathcal{E})\simeq 
Hom(a(\coprod X_{i}),\mathcal{E})$, pour toute famille finie de sch\'emas affines $X_{i}$). 
Il est clair que le morphisme $X \longrightarrow \mathcal{E}$ se factorise par 
$\mathcal{E}_{\nu}\subset \mathcal{E}$.
Enfin, il nous reste \`a remarquer que le morphisme induit
$$X \longrightarrow \mathcal{E}_{\nu}$$
est un \'epimorphisme Zariski local. En effet, pour tout anneau local $(A,m)$, et tout 
complexe parfait $E\in D(A)$ avec $Dim_{A/m}H^{i}(E\otimes_{A}^{\mathbb{L}}A/m) \leq \nu(i)$, 
il existe un complexe de $A$-modules libres $P^{*}$, avec 
$rang(P^{i})\leq \nu(i)$, et un quasi-isomorphisme $P^{*} \longrightarrow E$.
\hfill $\Box$ \\

On continue de se fixer un type $\nu$.
Pour $A\in k-CAlg$, on note $\mathcal{F}_{\nu}(A)$ le sous-ensemble de $dga(A)$ form\'ee
des $A$-dg-alg\`ebres propres et de type $\nu$. Cela d\`efinit un sous-pr\'efaisceau
$\mathcal{F}_{\nu}\subset dga$.

\begin{lem}\label{l4}
Le pr\'efaisceau $\mathcal{F}_{\nu}$ est quasi-compact.
\end{lem}

\textit{Preuve:} Nous allons appliquer le lemme \ref{l2} au morphisme naturel
$$\mathcal{F}_{\nu} \longrightarrow \mathcal{E}_{\nu}$$
qui \`a une dg-alg\`ebre ne retient que son complexe sous-jacent. Pour cela, 
on d\'efinit deux pr\'efaisceaux simpliciaux
$\underline{\mathcal{F}}_{\nu}$ et $\underline{\mathcal{E}}_{\nu}$
de la fa\c{c}on suivante. 

Pour $A\in k-CAlg$, on consid\`ere
la cat\'egorie $wA-dg-alg^{c}_{\nu}$, dont les objets sont 
les $A$-dg-alg\`ebres cofibrantes, propres de type $\nu$, et dont
les morphismes sont les \'equivalences de $A$-dg-alg\`ebres. Pour 
$A \mapsto A'$ un morphisme dans $k-CAlg$, on dispose d'un foncteur
de changement de bases
$$A'\otimes_{A} -  : wA-dg-alg^{c}_{\nu} \longrightarrow wA'-dg-alg^{c}_{\nu}.$$
Ceci permet de voir $A \mapsto wA-dg-alg^{c}_{\nu}$ comme un pseudo-foncteur
$k-CAlg \longrightarrow Cat$. A \'equivalence pr\`es, ce pseudo-foncteur se rectifie en un 
foncteur strict $k-CAlg \longrightarrow Cat$. En composant avec 
le foncteur nerf $N : Cat \longrightarrow SEns$ on obtient le pr\'efaisceau simplicial 
$\underline{\mathcal{F}}_{\nu}$. Par d\'efinition, $\underline{\mathcal{F}}_{\nu}(A)$ est naturellement
\'equivalent au nerf de la cat\'egorie $wA-dg-alg^{c}_{\nu}$, ce qui implique en particulier
qu'il existe un isomorphisme de pr\'efaisceaux $\pi_{0}(\underline{\mathcal{F}}_{\nu})\simeq \mathcal{F}_{\nu}$.

De m\^eme, pour $A\in k-CAlg$, on consid\`ere
la cat\'egorie $wC(A)^{c}_{\nu}$, dont les objets sont 
les complexes de $A$-modules cofibrants, parfaits de type $\nu$, et dont
les morphismes sont les quasi-isomorphismes. Pour 
$A \mapsto A'$ un morphisme dans $k-CAlg$, on dispose d'un foncteur
de changement de bases
$$A'\otimes_{A} -  : wC(A)^{c}_{\nu} \longrightarrow wC(A')^{c}_{\nu}.$$
Ceci permet de voir $A \mapsto wC(A)^{c}_{\nu}$ comme un pseudo-foncteur
$k-CAlg \longrightarrow Cat$. A \'equivalence pr\`es, ce pseudo-foncteur se rectifie en un 
foncteur strict $k-CAlg \longrightarrow Cat$. En composant avec 
le foncteur nerf $N : Cat \longrightarrow SEns$ on obtient le pr\'efaisceau simplicial 
$\underline{\mathcal{E}}_{\nu}$. Par d\'efinition, $\underline{\mathcal{E}}_{\nu}(A)$ est naturellement
\'equivalent au nerf de la cat\'egorie $wC(A)^{c}_{\nu}$, ce qui implique en particulier
qu'il existe un isomorphisme de pr\'efaisceaux $\pi_{0}(\underline{\mathcal{E}}_{\nu})\simeq \mathcal{E}_{\nu}$.

Enfin, il existe un morphisme de pr\'efaisceaux simpliciaux 
$$\underline{\mathcal{F}}_{\nu} \longrightarrow \underline{\mathcal{E}}_{\nu},$$
qui oublie la structure de dg-alg\`ebre. En appliquant les lemmes \ref{l2} et \ref{l3} on voit qu'il suffit
alors de montrer que pour tout sch\'ema affine $X=Spec\, A$, et tout 
morphisme $X \longrightarrow \underline{\mathcal{E}}_{\nu}$, le pr\'efaisceau
$\pi_{0}(\underline{\mathcal{F}}_{\nu}\times_{\underline{\mathcal{E}}_{\nu}}^{h}X)$ est quasi-compact.

On consid\`ere $\pi_{0}(\underline{\mathcal{F}}_{\nu}\times_{\underline{\mathcal{E}}_{\nu}}^{h}X)$ comme
un pr\'efaisceau sur $k-Aff/X\simeq A-CAlg$. D'apr\`es \cite[Thm. 1.2.15]{re}, ce pr\'efaisceau se d\'ecrit de la fa\c{c}on
suivante. Le morphisme $X \longrightarrow \underline{\mathcal{E}}_{\nu}$, correspond par 
le lemme de Yoneda \`a un complexes de $A$-modules $E$, cofibrant et de type $\nu$, que l'on 
pourra supposer \^etre strictement parfait (i.e. un complexe born\'e
de $A$-modules projectifs et de rangs finis). 
Alors, pour tout $A'\in A-CAlg$, on a un isomorphisme naturel
$$\pi_{0}(\underline{\mathcal{F}}_{\nu}\times_{\underline{\mathcal{E}}_{\nu}}^{h}X)(A')\simeq
\pi_{0}(Map_{Op}(C(k))(\mathcal{A}_{\infty},End(E)\otimes_{A}A'),$$
o\`u $Map_{Op(C(k))}$ d\'esigne l'espace des morphismes dans la cat\'egorie
des op\'erades au-dessus de la cat\'egorie monoidale $C(k)$, et o\`u 
$End(E)$ est l'op\'erade dans $C(A)$ des endomophismes de l'objet $E$. Or, comme l'op\'erade
$\mathcal{A}_{\infty}$ est cofibrante, et que toutes les op\'erades dans $C(k)$ sont fibrantes, 
le morphisme naturel
$$Hom_{Op(C(k))}(\mathcal{A}_{\infty},End(E)\otimes_{A}A') \longrightarrow 
\pi_{0}(Map_{Op(C(k))}(\mathcal{A}_{\infty},End(E)\otimes_{A}A'))$$
est surjectif. Pour terminer la preuve du lemme \ref{l4} il nous suffit de remarquer 
que le foncteur 
$A' \mapsto Hom_{Op(C(k))}(\mathcal{A}_{\infty},End(E)\otimes_{A}A')$
est repr\'esentable par un sch\'ema affine (ce qui se d\'eduit facilement du fait que
$E$ est un complexe born\'e de $A$-modules projectifs de rangs finis).
\hfill $\Box$ \\

\begin{lem}\label{l5}
Soit $A\in k-CAlg$, $B\in A-dg-alg$ une $A$-dg-alg\`ebre propre
et $M$ et $N$ deux $B$-modules parfaits sur $A$. 
Soit  $\mathcal{H}_{B}(M,N)$ le pr\'efaisceau sur $k-Aff/Spec\, A$ qui \`a $A'\in A-CAlg$ 
associe l'ensemble des morphismes  $M\otimes_{A}^{\mathbb{L}}A' \longrightarrow
N\otimes_{A}^{\mathbb{L}}A'$ dans la cat\'egorie $Ho(B\otimes_{A}^{\mathbb{L}}A'-Mod)$. Alors, 
$\mathcal{H}_{B}(M,N)$ est quasi-compact.
\end{lem}

\textit{Preuve:} On peut commencer par supposer que $B$ est une $A$-dg-alg\`ebre
cofibrante. De m\^eme, on supposera que 
$M$ est un $B$-module cofibrant. 

\begin{sublem}\label{sl}
Il existe un $B$-dg-module $P$ qui soit un complexe born\'e de $A$-modules projectifs
de rangs finis, et un isomorphisme dans $Ho(B-Mod)$,  $P \simeq N$.
\end{sublem}

\textit{Preuve du sous-lemme:}
Choissisons un complexe born\'e de $A$-modules projectifs libres et de rangs finis 
$P$ et un quasi-isomorphisme $u : P \longrightarrow N$. On factorise $u$ en 
$$\xymatrix{P \ar[r]^-{j} & N' \ar[r]^-{p} & N}$$
avec $j$ une cofibration triviale dans $C(A)$ et $p$ une fibration triviale dans $C(A)$. 
On consid\`ere $\underline{End}(p)$, la $A$-dg-alg\`ebre des endomophismes
du diagramme de complexes de $A$-modules $p : N' \rightarrow N$. Par d\'efinition, on a 
$$\underline{End}(p):=\underline{End}(N')\times_{\underline{Hom}(N',N)}\underline{End}(N).$$
On remarquera que le morphisme naturel 
$$\underline{End}(p) \longrightarrow \underline{End}(N)$$
est un une fibration triviale de $A$-dg-alg\`ebres. Comme $B$ est cofibrante, le morphisme
$B \longrightarrow \underline{End}(N)$ qui d\'etermine la structure de $B$-dg-module sur $N$
se rel\`eve \`a un morphisme $B \longrightarrow \underline{End}(p)$. En d'autres termes, il 
existe une structure de $B$-dg-module sur $N'$ tel que le morphisme $p : N' \longrightarrow N$
soit une \'equivalence. Ceci montre que l'on peut remplacer $N$ par $N'$ est 
donc supposer que le morphisme $u$ est une cofibration triviale de complexes de $A$-modules. 

Le morphisme
$$f : \underline{End}(u) \longrightarrow \underline{End}(N)$$
est alors une \'equivalences de $A$-dg-alg\`ebres. Ainsi, comme $B$ est cofibrante le morphisme
$\alpha : B \longrightarrow \underline{End}(N)$ se rel\`eve en un morphisme
$\beta : B \longrightarrow \underline{End}(u)$ \`a homotopie pr\`es. Soient $\Gamma^{1}(B)$ un objet
cylindre pour $B$ et 
$$h : \Gamma^{1}(B) \longrightarrow \underline{End}(N)$$
une homotopie entre $f\circ  \beta$ et $\alpha$. Le morphisme induit
$B \longrightarrow \underline{End}(u) \longrightarrow \underline{End}(P)$ d\'etermine
une structure de $B$-dg-module sur $P$. 
Le morphisme $h$ d\'etermine lui 
une structure de $\Gamma^{1}(B)$-dg-module sur le complexe $N$, et par les deux morphismes
naturels $B\rightrightarrows \Gamma^{1}(B)$, il d\'etermine aussi deux structures de $B$-dg-modules
sur le complexe $N$ (nous noterons $N_{1}$ et $N_{2}$ les deux $B$-dg-modules correspondants).
Le premi\`ere de ces structures, disons $N_{1}$, est la structure que l'on s'est initialement
donn\'ee sur $N$. La seconde, $N_{2}$,  est une structure de $B$-dg-modules tel que
le morphisme $u : P \longrightarrow N_{2}$ soit un morphisme de $B$-dg-modules. Enfin,
Comme les deux
morphismes $B\rightrightarrows \Gamma^{1}(B)$ sont des \'equivalences qui 
poss\`edent la projection naturelle $\Gamma^{1}(B) \longrightarrow B$ comme r\'etraction 
en commun, les deux objets $N_{1}$ et $N_{2}$ sont isomorphes dans $Ho(B-dg-Mod)$
(les deux foncteurs $Ho(\Gamma^{1}(B)) \rightrightarrows Ho(B-Mod)$ sont 
des quasi-inverses du m\^eme foncteur $-\otimes_{\Gamma^{1}(B)}^{\mathbb{L}}B$).
Ainsi, $P$ est finalement isomorphe \`a $N_{1}=N$ dans $Ho(B-dg-mod)$.
\hfill $\Box$ \\

On revient \`a la preuev du lemme \ref{l5}.
En utilisant le sous-lemme \ref{sl}, on supposera donc de plus que 
$N$ est un complexe born\'e de $A$-modules projectifs et de rangs finis. On consid\`ere alors
le pr\'efaisceau $Hom_{B}(M,N)$, qui \`a $A'\in A-CAlg$ associe
l'ensemble des morphismes de $B$-dg-modules
$M \longrightarrow N\otimes_{A}A'$. Comme le complexe
sous-jacent \`a $N$ est un complexe born\'e de $A$-modules
projectifs de rangs finis, le foncteur $Hom_{B}(M,N)$ est repr\'esentable
par un sch\'ema affine. Enfin, comme $M$ est cofibrant et que tout $B$-module est fibrant, 
le morphisme $Hom_{B}(M,N) \longrightarrow \mathcal{H}_{B}(M,N)$
est un \'epimorphisme de pr\'efaisceau. Ceci montre que $ \mathcal{H}_{B}(M,N)$
est quasi-compact.
\hfill $\Box$ \\

\begin{lem}\label{l7}
Soient $A\in k-CAlg$, $B\in A-dg-alg$ une $A$-dg-alg\`ebre propre, 
$M$ et $N$ deux $B$-modules parfaits sur $A$ et $f,g : M \longrightarrow N$ deux
morphismes de $B$-dg-modules. 
Soit  $\mathcal{EQ}(f,g)$ le sous-pr\'efaisceau de $Spec\, A$ 
qui \`a $A'\in k-CAlg$ 
associe le sous-ensemble de $X(A')$ form\'e des morphismes $A \longrightarrow A'$
tels que les deux morphismes
$$f\otimes_{A}^{\mathbb{L}}A', g\otimes_{A}^{\mathbb{L}}A' : M\otimes_{A}^{\mathbb{L}}A' \longrightarrow
N\otimes_{A}^{\mathbb{L}}A'$$ 
soient \'egaux dans la cat\'egorie $Ho(B\otimes_{A}^{\mathbb{L}}A'-Mod)$. Alors, 
$\mathcal{EQ}(f,g)$ est quasi-compact.
\end{lem}

\textit{Preuve:} La preuve suit le m\^eme principe que celle du lemme \ref{l5} et uitilise elle
aussi le sous-lemme \ref{sl}.

On peut supposer que $B$ est cofibrante, que $M$ est cofibrant, et de plus d'apr\`es le
sous-lemme \ref{sl} que $N$ est un complexe born\'e de $A$-modules projectifs et de rangs finis. 
Soit $\Gamma^{1}(M)$ un objet cyclindre pour le $B$-dg-module $M$. 
On consid\`ere alors $Hom_{B}(f,g)$, le pr\'efaisceau sur $k-Aff/Spec\, A$ qui \`a $A'\in A-CAlg$ associe
l'ensemble des morphismes $u : \Gamma^{1}(M) \longrightarrow N\otimes_{A}A'$, 
tels que les deux morphismes induits
$$M \rightrightarrows \Gamma^{1}(M) \longrightarrow N\otimes_{A}A'$$
soit \'egaux \`a $f\otimes_{A}A'$ et $g\otimes_{A}A'$. Le foncteur 
$Hom_{B}(f,g)$ est repr\'esentable par un sch\'ema affine au-dessus de $X$. De plus, 
le morphisme naturel
$$Hom_{B}(f,g) \longrightarrow X$$
a $\mathcal{EQ}(f,g)$ pour image. Ainsi, le pr\'efaisceau $\mathcal{EQ}(f,g)$ est quasi-compact.
\hfill $\Box$ \\

\begin{lem}\label{l8}
Soient $A\in k-CAlg$, $B\in A-dg-alg$ une $A$-dg-alg\`ebre propre et 
$d\geq 0$ un entier. Soit $X_{li,d}$ le sous-pr\'efaisceau de $X=Spec\, A$, 
qui \`a $A'\in k-CAlg$ associe le sous-ensemble de $X(A')$ form\'e
des morphismes $A\longrightarrow A'$ tels que 
$B\otimes_{A}^{\mathbb{L}}A'$ soit lisse et de dimension cohomologique
inf\'erieure \`a $d$ (relativement \`a $A'$). Alors, 
$X_{li,d}$ est quasi-compact.
\end{lem}

\textit{Preuve:} On consid\`ere $B$ comme un $B\otimes^{\mathbb{L}}_{A}B^{op}$-dg-module, et on lui 
applique la constuction de sa r\'esolution libre standard de \S 2.4 (on peut pour simplifier
supposer que $B$ est cofibrante).
On consid\`ere alors le morphisme de $B\otimes^{\mathbb{L}}_{A}B^{op}$-dg-modules
$$p : N:=|Sq_{d}\mathbb{R}_{*}(B)| \longrightarrow B.$$
Notons $F$ le pr\'efaisceau simplicial sur $k-Aff/Spec\, A$ qui \`a $A'\in A-CAlg$ associe
l'ensemble des morphismes $s : M \longrightarrow N\otimes_{A}^{\mathbb{L}}A'$
de $Ho(B\otimes^{\mathbb{L}}_{A}B^{op}-Mod)$ tels que $p\circ s=id$ 
(toujours dans $Ho(B\otimes^{\mathbb{L}}_{A}B^{op}-Mod)$). Le morphisme naturel
$F \longrightarrow X=Spec\, A$ a clairement pour pr\'efaisceau image $X_{li,d}$. De plus, les lemmes
\ref{l5} et \ref{l7} impliquent que $F$ est quasi-compact. Ainsi, $X_{li,d}$ est quasi-compact.
\hfill $\Box$ \\

\textit{Preuve du th\'eor\`eme \ref{t2}:} Notons $\mathcal{F}_{\nu,d}$ le pr\'efaisceau sur $k-Aff$
qui \`a $A$ associe l'ensemble des classes d'isomorphismes
de $A$-dg-alg\`ebres $B$ dans $Ho(A-dg-alg)$ qui sont propres et lisses et qui 
v\'erifient les deux conditions suivantes
\begin{enumerate}
\item $B$ est de type $\nu$ (relativement \`a $A$). 
\item Il existe $a_{1}, \dots, a_{n}\in A$ avec 
$\sum a_{i}=1$, tels que chaque $B\otimes_{A}A[a_{i}^{-1}]$ soit de dimension cohomologique 
inf\'erieure \`a $d$ (relativement \`a $A[a_{i}^{-1}]$).
\end{enumerate}

Il s'agit bien entendu de montrer que le pr\'efaisceau $\mathcal{F}_{\nu,d}$
est quasi-compact. 

Le pr\'efaisceau $\mathcal{F}_{\nu,d}$ est un sous-pr\'efaisceau 
de $\mathcal{F}_{\nu}$ consid\'er\'e dans le lemme \ref{l4}. Soit $X=Spec\, A \longrightarrow \mathcal{F}_{\nu}$
un \'epimorphisme Zariski local. Ce morphisme correspond \`a une
$A$-dg-alg\`ebre $B$ propre et de type $\nu$. Par le lemme \ref{l8} $X_{li,d}$, le lieu o\`u $B$ 
est lisse de dimension cohomologique inf\'erieure \`a $d$ est quasi-compact. De toute
\'evidence, le morphisme $X_{li,d} \longrightarrow \mathcal{F}_{\nu}$ se factorise
par $\mathcal{F}_{\nu,d}$, et le morphisme induit $X_{li,d} \longrightarrow \mathcal{F}_{\nu,d}$
est un \'epimorphisme Zariski local.
Ceci termine la preuve du th\'eor\`eme. 
\hfill $\Box$ \\

\section{Le th\'eor\`eme de finitude homotopique}

Nous arrivons enfin au th\'eo\`eme de finitude homotopique qui s'\'ennonce comme suit.

\begin{thm} \label{t3}
Soit $k$ un anneau commutatif. 
Soient $\nu$ un type et 
$d\in \mathbb{N}$. Alors, il existe un entier $n(\nu,d)$ qui v\'erifie les deux propri\'et\'es suivantes.
\begin{enumerate}
\item Pour toute $k$-dg-alg\`ebre $B$ propre et lisse, de type $\nu$ et de dimension cohomologique inf\'erieure \`a $d$, 
il existe des \'el\'ements $f_{1}, \dots, f_{n}\in k$ avec $\sum_{i}f_{i}=1$ et tels que
pour tout $i$ les morphismes naturels
$$\mathbb{L}(i_{n(\nu,d)})_{!}i_{n(\nu,d)}^{*}(B\otimes_{k}k[f_{i}^{-1}]) \longrightarrow B\otimes_{k}k[f_{i}^{-1}]$$
poss\`ede une section dans $Ho(k[f_{i}^{-1}]-dg-alg)$. 
\item Pour deux $k$-dg-alg\`ebres $B$ et $B'$ propres et lisses, de type $\nu$ et de dimension 
cohomologique inf\'erieure \`a $d$, il existe des \'el\'ements $f_{1}, \dots, f_{n}\in k$ avec $\sum_{i}f_{i}=1$ et tels que
pour tout $i$
les morphismes d'ensembles simpliciaux
$$Map_{k[f_{i}^{-1}]-dg-alg}(B\otimes_{k}k[f_{i}^{-1}],B'\otimes_{k}k[f_{i}^{-1}]) \longrightarrow Map_{k[f_{i}^{-1}]-\mathcal{A}_{n(\nu,d)}-alg}(i_{n(\nu,d)}^{*}(B\otimes_{k}k[f_{i}^{-1}]),i_{n(\nu,d)}^{*}(B'\otimes_{k}k[f_{i}^{-1}]))$$
$$Map^{eq}_{k[f_{i}^{-1}]-dg-alg}(B\otimes_{k}k[f_{i}^{-1}],B'\otimes_{k}k[f_{i}^{-1}]) \longrightarrow Map^{eq}_{k[f_{i}^{-1}]-\mathcal{A}_{n(\nu,d)}-alg}(i_{n(\nu,d)}^{*}(B\otimes_{k}k[f_{i}^{-1}]),i_{n(\nu,d)}^{*}(B'\otimes_{k}k[f_{i}^{-1}]))$$
poss\`edent des r\'etractions dans la cat\'egorie $Ho(SEns)$. 
\end{enumerate}
\end{thm}

Avant d'attaquer la preuve du th\'eor\`eme nous en d\'eduisons le corollaire suivant.

\begin{cor}\label{c3}
Soit $k$ un anneau local commutatif. 
Pour tout type $\nu$ e tout entier $d\geq 0$, il existe un entier $n(\nu,d)$ 
qui v\'erifie les conditions suivantes. 
\begin{enumerate}
\item Pour toute $k$-dg-alg\`ebre $B$ propre et lisse, de type $\nu$ et de dimension cohomologique inf\'erieure \`a $d$, 
les morphismes naturels
$$\mathbb{L}(i_{n(\nu,d)})_{!}i_{n(\nu,d)}^{*}(B) \longrightarrow B$$
poss\`ede une section dans $Ho(k-dg-alg)$. 
\item Pour deux $k$-dg-alg\`ebres $B$ et $B'$ propres et lisses, de type $\nu$ et de dimension 
cohomologique inf\'erieure \`a $d$, 
les morphismes d'ensembles simpliciaux
$$Map_{k-dg-alg}(B,B') \longrightarrow Map_{k-\mathcal{A}_{n(\nu,d)}-alg}(i_{n(\nu,d)}^{*}(B),i_{n(\nu,d)}^{*}(B'))$$
$$Map^{eq}_{k-dg-alg}(B,B') \longrightarrow Map^{eq}_{k-\mathcal{A}_{n(\nu,d)}-alg}(i_{n(\nu,d)}^{*}(B),i_{n(\nu,d)}^{*}(B'))$$
poss\`edent des r\'etractions dans la cat\'egorie $Ho(SEns)$. 

\item 
Deux $k$-dg-alg\`ebres $B$ et $B'$ 
de type $\nu$ et de dimension cohomologique inf\'erieure \`a $d$ sont 
isomorphes dans $Ho(k-dg-alg)$ si et seulement si 
les $\mathcal{A}_{n(\nu,d)}$-alg\`ebres $i_{n(\nu,d)}^{*}(B)$ et $i_{n(\nu,d)}^{*}(B')$ sont isomorphes
dans $Ho(\mathcal{A}_{n(\nu,d)}-dg-alg)$. 
\end{enumerate}
\end{cor}

\subsection{Petit d\'etour par la ''nouvelle g\'eom\'etrie alg\'ebrique courageuse''}

Pour ce paragraphe on se fixe un anneau commutatif $k$. \\

Soit $Hk-Mod$ la cat\'egorie de mod\`eles des $Hk$-modules dans
les spectres sym\'etriques au sens de \cite{hss}. On munira 
$Hk-Mod$ de sa structure de cat\'egorie de mod\`eles stbale et positive
de \cite{s}. Le smash produit de spectres sym\'etriques induit une structure
monoidale sym\'etrique $-\wedge_{Hk} -$ sur $Hk-Mod$ qui en fait une
cat\'egorie de mod\`ele monoidale v\'erifiant l'axiome du monoide. Enfin, on notera 
$Hk-CAlg$ la cat\'egorie des monoides commutatifs dans $Hk-Mod$ (i.e. 
des $Hk$-alg\`ebres commutatives). D'apr\`es \cite{s} elle est munie d'une structure
de cat\'egorie de mod\`eles pour laquelle les fibrations et les \'equivalences
sont d\'efinies dans $Hk-Mod$. De m\^eme, nous noterons $Hk-Alg$ la cat\'egorie
des monoides dans $Hk-Mod$ (i.e. des $Hk$-alg\`ebres). Elle est elle aussi munie
de sa structure de mod\`eles  pour laquelle les fibrations et les \'equivalences
sont d\'efinies dans $Hk-Mod$. On sait qu'il existe une chaine
d'\'equivalences de Quillen entre $Hk-Alg$ et $k-dg-alg$ (voir \cite{s2}). 
On identifiera les cat\'egories homotopiques $Ho(Hk-Alg)$ et 
$Ho(k-dg-alg)$ \`a travers cette \'equivalence, ainsi que les
ensembles simpliciaux de morphismes correspondant. 

Notons $Hk-Aff$ la cat\'egorie oppos\'ee \`a $Hk-CAlg$, et
$SPr(Hk-Aff)$ la cat\'egorie des foncteurs $Hk-CAlg \longrightarrow SEns$. 
Nous munirons $SPr(Hk-Aff)$ de sa structure de mod\`ele projective niveaux 
par niveaux pour laquelle les \'equivalences et les fibrations sont d\'efinies
niveaux par niveaux au-dessus des objets de $Hk-CAlg$. 
Pour une $Hk$-alg\`ebre commutative $A\in Hk-CAlg$, on pose 
$$\mathbb{R}Spec\, A : Hk-CAlg \longrightarrow SEns,$$
 qui \`a $B\in Hk-CAlg$ fait correspondre l'ensemble simplicial des morphismes
de $Q(A)$ dans $R(B)$
($Hk-CAlg$ est une cat\'egorie de mod\`eles simpliciale)
$$(\mathbb{R}Spec\, A):=\underline{Hom}(Q(A),R(B)),$$
o\`u $Q$ et $R$ sont des foncteurs de remplacement cofibrant
et fibrant dans $Hk-CAlg$. Ceci d\'efinit un foncteur au niveau des cat\'egories
homotopiques
$$\mathbb{R}Spec : Ho(Hk-CAlg)^{op} \longrightarrow Ho(SPr(Hk-Aff)).$$
Le lemme de Yoneda pour la cat\'egories de mod\`eles (voir e.g. \cite{hagI}) implique que
$\mathbb{R}Spec$ est un foncteur pleinement fid\`ele. 

Un objet dans l'image essentielle du foncteur 
$\mathbb{R}Spec$ sera dit \emph{affine}. Comme le
foncteur $\mathbb{R}Spec$ commute avec les limites homotopiques et qu'il
est pleinement fid\`ele, on voit qu'un foncteur $Hk-CAlg \longrightarrow SEns$
qui est limite homotopique d'affines est lui-m\^eme affine. \\

\begin{df}\label{d8}
Un foncteur $F : Hk-CAlg \longrightarrow SEns$
est \emph{homotopiquement de pr\'esentation finie} si 
pour tout syst\`eme filtrant d'objets $\{A_{\alpha}\}$ dans
$Hk-CAlg$, le morphisme naturel
$$Hocolim_{\alpha} F(A_{\alpha})\simeq  Colim_{\alpha} F(A_{\alpha}) \longrightarrow F(Hocolim_{\alpha}A_{\alpha})$$
est une \'equivalence. 
\end{df}

Soit maintenant $B$ et $B'$ deux $k$-dg-alg\`ebres propres. On les voit \`a travers l'\'equivalence
$Ho(k-dg-alg)\simeq Ho(Hk-Alg)$ comme des $Hk$-alg\`ebres. On d\'efinit alors
un foncteur
$$Eq(B,B') : Hk-CAlg \longrightarrow SEns,$$
qui \`a $A\in Hk-CAlg$ associe 
$$Eq(B,B')(A):=\underline{Hom}^{eq}_{Hk-Alg}(Q(B)\wedge_{Hk}A,R(Q(B')\wedge_{Hk}A)),$$
o\`u $Q$ et $R$ sont des foncteurs de remplacement cofibrant et fibrant dans
$Hk-Alg$. Ce foncteur sera consid\'er\'e comme un objet dans $Ho(SPr(Hk-Aff))$. 

\begin{lem}\label{l9}
Le foncteur $Eq(B,B')$ est affine. 
\end{lem}

\textit{Preuve:} La preuve de cette proposition se fait en deux \'etapes. On commence
par consid\'erer un foncteur auxiliaire
$$Hom(B,B') : Hk-CAlg \longrightarrow SEns,$$
qui \`a $A\in Hk-CAlg$ associe 
$$Hom(B,B')(A):=\underline{Hom}_{Hk-Alg}(Q(B),R(Q(B')\wedge_{Hk}A))=
\underline{Hom}_{HA-Alg}(Q(B)\wedge_{Hk}A,R(Q(B')\wedge_{Hk}A)).$$
Ainsi, le foncteur $Eq(B,B')$ est un sous-foncteur de $Hom(B,B')$. 

Commen\c{c}ons pas montrer que $Hom(B,B')$ est affine. 
Pour cela, on \'ecrit $B$ comme une colimite $Colim_{n}B_{n}$, o\`u 
pour tout entier $n$ il existe un carr\'e cocart\'esien dans $k-dg-alg$
$$\xymatrix{
B_{n} \ar[r] & B_{n+1} \\
\coprod_{i} S_{k}(n_{i}) \ar[u] \ar[r] & \coprod_{i}D_{k}(n_{i}+1), \ar[u]}$$
o\`u le coproduit est pris sur l'ensembles des paires $(n_{i},f)$, avec
$n_{i}\in \mathbb{Z}$ et $f : S_{k}(n_{i}) \longrightarrow B_{n}$
est un morphisme (voir \cite[\S 2]{ho}). Dans ce cas, le foncteur 
$Hom(B,B')$ est \'equivalent \`a la limite homotopique 
$Holim_{n}Hom(B_{n},B')$. Pour tout $n$ on dispose de plus d'un diagrmme homotopiquement cart\'esien
$$\xymatrix{
Hom(B_{n+1},B') \ar[r] \ar[d] & Hom(B_{n+1},B')\ar[d] \\
\prod_{i} Hom(D_{k}(n_{i}+1),B')  \ar[r] & \prod_{i}Hom(S_{k}(n_{i}),B').}$$
Comme les objets affines sont stable par limites homotopiques il 
nous suffit de v\'erifier que pour un entier $n$ les
objets $Hom(D_{k}(n),B') $ et $Hom(S_{k}(n),B')$ sont affines. 
Comme $D_{k}(n)$ est \'equivalente \`a $k$, on a
$Hom(D_{k}(n),B')\simeq \mathbb{R}Spec\, Hk$. 
La $k$-dg-alg\`ebre $S_{k}(n)$ est libre engendr\'e par 
un \'el\'ement en degr\'e $-n$, on a pour tout $A\in Hk-CAlg$ 
$$Hom(S_{k}(n),B')(A)\simeq Map_{Hk-Mod}(Hk,B'\wedge_{Hk}^{\mathbb{L}}A[-n]).$$
Notons $(B')^{\vee}:=\mathbb{R}\underline{Hom}(B',Hk)$ le $Hk$-module
dual (d\'eriv\'e) de $B'$ (notons que comme $B'$ est propre, il est fortement dualisable
comme $Hk$-module,  voir \cite{hagII} pour la notion de fortement dualizable). Consid\'erons 
$A_{0}$ la $Hk$-alg\`ebre commutative libre engendr\'ee 
par $(B')^{\vee}[n]$. Par dualit\'e, on a pour $A\in Hk-CAlg$
$$(\mathbb{R}Spec\, A_{0})(A)\simeq Map_{Hk-Mod}((B')^{\vee}[n],A)\simeq
Map_{Hk-Mod}(Hk,B'\wedge_{Hk}^{\mathbb{L}}A[-n]).$$
Le foncteur $Hom(S_{k}(n),B')$ est donc isomorphe dans $Ho(SPr(Hk-Aff))$ 
\`a $\mathbb{R}Spec\, A_{0}$. 

L'\'equivalence $\mathbb{R}Spec\ A_{0} \simeq Hom(B,B')$
d\'efinit un morphisme de $A_{0}$-alg\`ebres 
$u : Q(B)\wedge_{Hk}A_{0} \longrightarrow R(Q(B')\wedge_{Hk}A_{0})$. Notons 
$K$ la fibre homotopique de ce morphisme dans $Ho(A_{0}-Mod)$. 
Le sous-foncteur 
$$Eq(B,B') \subset Hom(B,B')\simeq \mathbb{R}Spec\, A_{0}$$ 
est tel que pour tout $A\in Hk-CAlg$, $Eq(B,B')(A_{0})$  est le sous-ensemble
simplicial de $Map_{Hk-CAlg}(A_{0},A)$ form\'e des morphismes
$v : A \longrightarrow A_{0}$ tels que $K\wedge_{A}^{\mathbb{K}}A_{0}\simeq *$. 
Comme les $Hk$-modules $B$ et $B'$ sont parfaits, le sous-foncteur
$Eq(B,B')$ est affine de la forme $\mathbb{R}Spec\, (A_{0})_{K}$ (voir
\cite[\S 1.2.10]{hagII}).  \hfill $\Box$ \\

\begin{lem}\label{l10}
Si $B$ et $B'$ sont propres et lisses alors
le foncteur $Eq(B,B')$ est homotopiquement de pr\'esentation finie. 
\end{lem}

\textit{Preuve:} Soit $\{A_{\alpha}\}$ un syst\`eme filtrant 
dans $Hk-CAlg$ de colimite homotopique $A=Hocolim_{\alpha}$. On consid\`ere le carr\'e commutatif
$$\xymatrix{
Colim_{\alpha}Eq(B,B')(A_{\alpha}) \ar[r] \ar[d] & Eq(B,B')(A) \ar[d] \\
Colim_{\alpha}Hom(B,B')(A_{\alpha}) \ar[r] & Hom(B,B')(A).}$$
Par d\'efinition, les morphismes verticaux induisent des isomorphismes
sur les $\pi_{i}$ pour $i>0$ et des injections sur les $\pi_{0}$. 
De plus, $B$ \'etant homotopiquement de pr\'esentation finie on a 
$$
Colim_{\alpha}Hom(B,B')(A_{\alpha}) \simeq Colim_{\alpha}Map_{Hk-Alg}(B,B'\wedge_{A}^{\mathbb{L}}A_{\alpha})
\simeq Map_{Hk-Alg}(B,Hocolim_{\alpha} B'\wedge_{A}^{\mathbb{L}}A_{\alpha})$$
$$\simeq Map_{Hk-Alg}(B,B'\wedge_{A}^{\mathbb{L}} (Hocolim_{\alpha} A_{\alpha}))\simeq
Map_{Hk-Alg}(B,B'\wedge_{A}^{\mathbb{L}}A)\simeq 
Hom(B,B')(A).$$
Ainsi, on en d\'eduit que le morphisme horizontal du haut 
induit des isomorphismes sur les $\pi_{i}$ pour $i>0$ et
une injection sur $\pi_{0}$. Il nous reste \`a voir que 
ce morphisme est aussi surjectif sur les $\pi_{0}$. Soit 
$u : B\wedge_{Hk}^{\mathbb{L}}A \simeq B\wedge_{Hk}^{\mathbb{L}}A$
un isomorphisme dans $Ho(A-Alg)$. Notons $v$ un inverse 
de $u$ dans $Ho(A-Alg)$. Comme 
$$[B'\wedge_{Hk}^{\mathbb{L}}A,B\wedge_{Hk}^{\mathbb{L}}A]\simeq 
Colim_{\alpha}[B'\wedge_{Hk}^{\mathbb{L}}A_{\alpha},B\wedge_{Hk}^{\mathbb{L}}A_{\alpha}]$$
il existe un $\alpha$ tel que $u$ et $v$ soit d\'efini sur $A_{\alpha}$ (nous
noterons ces morphismes $u_{\alpha}$ et $v_{\alpha}$ dans $Ho(A_{\alpha}-Alg)$).
De plus, on a $uv=id$ dans $[B'\wedge_{Hk}^{\mathbb{L}}A,B'\wedge_{Hk}^{\mathbb{L}}A]$. 
Or, comme nous l'avons vu ci-dessus on a 
$$[B'\wedge_{Hk}^{\mathbb{L}}A,B'\wedge_{Hk}^{\mathbb{L}}A]\simeq 
Colim_{\alpha}[B'\wedge_{Hk}^{\mathbb{L}}A_{\alpha},B'\wedge_{Hk}^{\mathbb{L}}A_{\alpha}].$$
Ainsi, on peut choissir $\alpha$ de sorte \`a ce que $u$ et $v$ soit
d\'efinis sur $A_{\alpha}$, et de plus $u_{\alpha}v_{\alpha}=id$. 
De m\^eme, on trouve que l'on peut choisir $\alpha$ de sorte \`a ce que
de plus $v_{\alpha}u_{\alpha}=id$. Ainsi, $u_{\alpha}$ est un isomorphisme
dans $Ho(A_{\alpha}-Alg)$. Ceci termine la preuve du lemme.
\hfill $\Box$ \\

\begin{cor}\label{c4}
Soient $B$ et $B'$ deux $k$-dg-alg\`ebres propres et lisses. Alors il existe 
un entier $n$ tel que pour tout $A\in k-CAlg$ le morphisme d'ensembles simpliciaux
$$Map^{eq}_{A-dg-alg}(B\otimes_{k}^{\mathbb{L}}A,B'\otimes_{k}^{\mathbb{L}}A) \longrightarrow 
Map^{eq}_{k-\mathcal{A}_{n}-alg}(i_{n}^{*}(B)\otimes_{k}^{\mathbb{L}}A,i_{n}^{*}(B')\otimes_{k}^{\mathbb{L}}A)$$
poss\`ede une r\'etaction dans $Ho(SEns)$. 
\end{cor}

\textit{Preuve:} Nous avons vu que le foncteur $Eq(B,B')$ est affine et homotopiquement de
pr\'esentation finie (voir Lem. \ref{l9} et \ref{l10}). De m\^eme, on peut d\'efinir un foncteur $Eq_{n}(B,B')$, qui 
\`a $A\in Hk-Alg$ associe 
$\underline{Hom}_{Hk-\mathcal{A}_{n}-alg}(Q(i_{n}^{*}(B)),R(Q(i_{n}^{*}(B'))\wedge_{Hk}A))$.
On laisse le soin aux lecteurs de g\'en\'eraliser les constructions pr\'ec\'edentes
des alg\`ebres associatives aux alg\`ebres sur une op\'erade dans $Hk-Mod$. 
En particulier, de fa\c{c}on tout \`a fait analogue au lemme \ref{l9} on montre que 
$Eq_{n}(B,B')$ est affine. 
Notons $A_{0} \in Hk-CAlg$ tel que $\mathbb{R}Spec\, A_{0}\simeq Eq(B,B')$. De m\^eme, 
notons $A_{n} \in Hk-CAlg$ tel que $\mathbb{R}Spec\, A_{n}\simeq Eq_{n}(B,B')$. 
D'apr\`es la proposition \ref{p1} on a
$$Eq(B,B')\simeq Holim_{n}Eq_{n}(B,B'),$$
et donc, comme le foncteur $\mathbb{R}Spec$ est pleinement fid\`ele, on a
$$A_{0}\simeq Hocolim_{n}A_{n}.$$
Enfin, comme $Eq(B,B')$ est homotopiquement de pr\'esentation finie, $A_{0}$ est 
homotopiquement de pr\'esentation finie dans $Hk-CAlg$, et donc il existe un entier 
$n$ tel que le morphisme naturel 
$A_{n} \longrightarrow A$ poss\`ede une section dans $Ho(Hk-CAlg)$. Ainsi, 
$$Eq(B,B')\simeq \mathbb{R}Spec\, A \longrightarrow Eq_{n}(B,B')\simeq \mathbb{R}Spec\, A_{n}$$
poss\`ede une r\'etraction dans $Ho(SPr(Hk-Aff))$. En prenant les valeurs
sur $A\in Hk-CAlg$ on trouve que le r\'esultat annonc\'e.
\hfill $\Box$ \\

\subsection{Preuve du th\'eor\`eme}

Nous sommes maintenant en mesure de d\'emontrer notre th\'eor\`eme \ref{t3}. \\

Fixons un type $\nu$ et un entier $d\geq 0$. D'apr\`es le th\'eor\`eme
\ref{t2} il existe une $k$-alg\`ebre commutative $A_{0}$ et une
$A$-dg-alg\`ebre $B_{0}$ propre de type $\nu$, lisse et de dimension cohomologique inf\'erieure \`a $d$
qui v\'erifie la conditions suivante: pour toute $k$-dg-alg\`ebre $B$ propre, lisse et 
de dimension cohomologique inf\'erieure \`a $d$, il existe des \'el\'ements 
$f_{1}, \dots f_{m} \in k$, et des morphismes $A_{0} \longrightarrow k[f_{i}^{-1}]$ tel que
$B_{0}\otimes^{\mathbb{L}}_{A_{0}}k[f_{i}^{-1}]$ et $B\otimes_{k}k[f_{i}^{-1}]$ 
soient isomorphes dans $Ho(k[f_{i}^{-1}]-dg-alg)$ (pour tout $i$). 

On consid\`ere $A_{1}:=A_{0}\otimes_{k}A_{0}$, et deux $A_{1}$-dg-alg\`ebres
$$B_{1}:=B_{0}\otimes_{A_{0}}^{\mathbb{L}}(A_{0}\otimes_{k}A_{0}) \qquad 
B_{2}:=(A_{0}\otimes_{k}A_{0})\otimes_{A_{0}}^{\mathbb{L}}B.$$

D'apr\`es le corollaire \ref{c2}, et le corollaire \ref{c4} on sait qu'il existe un entier $n(\nu,d)$ qui v\'erifie les deux conditions suivantes.
\begin{enumerate}
\item Le morphisme
$$\mathbb{L}(i_{n})_{!}i_{n}^{*}(B_{0}) \longrightarrow B_{0}$$
poss\`ede une section dans $Ho(A_{0}-dg-alg)$.
\item Pour tout $A\in A_{1}-CAlg$, le morphisme 
$$Map^{eq}_{A-dg-alg}(B_{1}\otimes_{A_{1}}^{\mathbb{L}}A,B_{2}\otimes_{A_{1}}^{\mathbb{L}}A) \longrightarrow 
Map^{eq}_{A-\mathcal{A}_{n}-alg}(i_{n}^{*}(B_{1})\otimes_{A_{1}}^{\mathbb{L}}A,
i_{n}^{*}(B_{2})\otimes_{A_{1}}^{\mathbb{L}}A)$$
poss\`ede une r\'etaction dans $Ho(SEns)$. 
\end{enumerate}

La propri\'et\'e semi-universelle du couple $(A_{0},B_{0})$ rappel\'ee plus haut implique alors le point 
$(1)$ et la seconde partie du point $(2)$ du th\'eor\`eme (celle concernant les
$Map^{eq}$). Cela dit, le point $(1)$ du th\'eor\`eme implique aussi 
la premi\`ere partie du point $(2)$ \`a l'aide de l'\'equivalence
$$Map_{k-\mathcal{A}_{n}-alg}(i_{n}^{*}(B),i_{n}^{*}(B'))\simeq
Map_{k-dg-alg}(\mathbb{L}(i_{n})_{!}i^{*}_{n}(B),B').$$
Ceci termine la preuve du th\'eor\`eme \ref{t3}. 
\hfill $\Box$ \\

En corollaire de la preuve pr\'ec\'edente on voit que 
l'entier $n(\nu,d)$ peut \^etre choisi ind\'ependant du l'anneau $k$. 
En effet, si $X=Spec\, A \longrightarrow \mathcal{F}_{\nu,d,\mathbb{Z}}$
un \'epimorphisme Zariksi local, o\`u $\mathcal{F}_{\nu,d,\mathbb{Z}} \in \widehat{\mathbb{Z}-Aff}$ est 
d\'efini comme dans la preuve du th\'eor\`eme \ref{t2} mais pour $k=\mathbb{Z}$. Alors, 
le morphisme induit $X \times Spec\, k \longrightarrow \mathcal{F}_{\nu,d,\mathbb{Z}}\times Spec\, k$ 
est encore un \'epimorphisme Zariski local dans $\widehat{k-Aff}\simeq \widehat{\mathbb{Z}-Aff/Spec\, k}$. 
Or, $\mathcal{F}_{\nu,d,\mathbb{Z}}\times Spec\, k$ est le pr\'efaisceau d\'efini dans la preuve du th\'eor\`eme \ref{t2}. 
Ainsi, le couple $(A_{0},B_{0})$ utilis\'e en d\'ebut de preuve du th\'eor\`eme
\ref{t3} ci-dessus peut \^etre choisi de la forme $(A\otimes_{\mathbb{Z}} k,B\otimes_{A}^{\mathbb{L}}(A\otimes_{\mathbb{Z}} k))$, 
pour $B$ une certaine $A$-dg-alg\`ebre propre et lisse. La version finale du th\'eor\`eme
de finitude homotopique est alors la suivante.

\begin{thm} \label{t4}
Soient $\nu$ un type et 
$d\in \mathbb{N}$. Alors, il existe un entier $n(\nu,d)$ tel que pour tout anneau commutatif $k$ 
les deux propri\'et\'es suivantes sont satisfaites.
\begin{enumerate}
\item Pour toute $k$-dg-alg\`ebre $B$ propre et lisse, de type $\nu$ et de dimension cohomologique inf\'erieure \`a $d$, 
il existe des \'el\'ements $f_{1}, \dots, f_{n}\in k$ avec $\sum_{i}f_{i}=1$ et tels que
pour tout $i$ les morphismes naturels
$$\mathbb{L}(i_{n(\nu,d)})_{!}i_{n(\nu,d)}^{*}(B\otimes_{k}k[f_{i}^{-1}]) \longrightarrow B\otimes_{k}k[f_{i}^{-1}]$$
poss\`ede une section dans $Ho(k[f_{i}^{-1}]-dg-alg)$. 
\item Pour deux $k$-dg-alg\`ebres $B$ et $B'$ propres et lisses, de type $\nu$ et de dimension 
cohomologique inf\'erieure \`a $d$, il existe des \'el\'ements $f_{1}, \dots, f_{n}\in k$ avec $\sum_{i}f_{i}=1$ et tels que
pour tout $i$
les morphismes d'ensembles simpliciaux
$$Map_{k[f_{i}^{-1}]-dg-alg}(B\otimes_{k}k[f_{i}^{-1}],B'\otimes_{k}k[f_{i}^{-1}]) \longrightarrow Map_{k[f_{i}^{-1}]-\mathcal{A}_{n(\nu,d)}-alg}(i_{n(\nu,d)}^{*}(B\otimes_{k}k[f_{i}^{-1}]),i_{n(\nu,d)}^{*}(B'\otimes_{k}k[f_{i}^{-1}]))$$
$$Map^{eq}_{k[f_{i}^{-1}]-dg-alg}(B\otimes_{k}k[f_{i}^{-1}],B'\otimes_{k}k[f_{i}^{-1}]) \longrightarrow Map^{eq}_{k[f_{i}^{-1}]-\mathcal{A}_{n(\nu,d)}-alg}(i_{n(\nu,d)}^{*}(B\otimes_{k}k[f_{i}^{-1}]),i_{n(\nu,d)}^{*}(B'\otimes_{k}k[f_{i}^{-1}]))$$
poss\`edent des r\'etractions dans la cat\'egorie $Ho(SEns)$. 
\end{enumerate}
\end{thm}

\end{document}